\documentclass[10pt]{article}

\usepackage{amsmath}
\usepackage{amssymb}
\usepackage{indentfirst}
\usepackage{graphics} 
\usepackage{color}
\makeatletter

\@addtoreset{equation}{section} 
\makeatother
\def\<{\langle}
\def\>{\rangle}

\newtheorem{lem}{Lemma}[section]
\newtheorem{theo}{Theorem}[section]
\newtheorem{rem}{Remark}[section]

\makeatletter
   
   \@addtoreset{equation}{section}
\makeatother

\begin{document}
\title{\bf  Strongly damped wave equations with  mass-like\\ terms of the logarithmic-Laplacian}

\author{Alessandra Piske\thanks {alessandrapiske@gmail.com}\, \, and \,Ruy Coimbra Char\~ao\thanks{ ruy.charao@ufsc.br}
\\{ \small Graduate Program in Pure and Applied Mathematics}  
 \\{\small Department of Mathematics} \\{\small Federal University of Santa Catarina} \\ {\small 88040-900, Florianopolis, Brazil,} 
\\
and\\Ryo Ikehata\thanks{Corresponding author: ikehatar@hiroshima-u.ac.jp} \\ {\small Department of Mathematics, Division of Educational Sciences}\\ {\small Graduate School of Humanities and Social Sciences} \\ {\small Hiroshima University} \\ {\small Higashi-Hiroshima 739-8524, Japan}}
\date{}
\maketitle
%%%%%%%%%%%%%%%%%%%%%%%%%E½@abstractE½@%%%%%%%%%%%%%%%%%%%%%%
\begin{abstract}
We consider strongly damped wave equations with logarithmic mass-like terms with a parameter $\theta \in (0,1]$. This research is a part of a series of wave equations that was initiated by Char\~ao-Ikehata \cite{Log-damping},  Char\~ao-D'Abbicco-Ikehata considered in \cite{CDI} depending on a parameter $\theta \in (1/2,1)$ and Piske-Char\~ao-Ikehata \cite{JDE-2022} for small parameter $\theta \in (0, 1/2)$. We derive a leading term as $t \to \infty$ of the solution, and by using it, a growth and a decay property of the solution itself can be precisely studied in terms of $L^{2}$-norm. An interesting aspect appears in the case of $n = 1$, roughly speaking, a small $\theta$ produces a diffusive property, and a large $\theta$ gives a kind of singularity, expressed by growth rates. 
%We first consider  strongly damped wave equations with logarithmic dispersion term with a parameter $\theta \in (0,1]$. This research is a part of a series of wave equations that was initiated by Char\~ao-Ikehata \cite{Log-damping},  Char\~ao-D'Abbicco-Ikehata considered in \cite{CDI} depending on a parameter $\theta \in (1/2,1)$ and Piske-Char\~ao-Ikehata \cite{JDE-2022} for small parameter $\theta \in (0, 1/2)$. We derive the leading term as $t \to \infty$ of the solution, and by using it, an infinite time blowup and a decay property of the solution itself can be precisely studied in terms of $L^{2}$-norm. An interesting aspect appears in the case of $n = 1$, roughly speaking, a small $\theta$ produces a diffusive property, and a large $\theta$ gives a kind of singularity, expressed by growth rates. 
\end{abstract}
%%%%%%%%%%%%%%%%%%E½@introductionE½@%%%%%%%%%%%%%%%%%%%%%%%
\vspace{-1cm}

\section{Introduction}
We consider in this work strongly damped wave equations under effects of a logarithmic mass-like term dependng on a parameter $\theta$:
%%%
\begin{align}
&u_{tt} - \Delta u +m^2 L_\theta u - \Delta u_t =0, \quad t>0, \quad x \in {\bf R}^n, \label{equation}  \\
&u(0,x) = u_0(x), \quad u_t(0,x)= u_1(x), \quad x \in {\bf R}^n, \label{initial}
\end{align}
where $\theta \in (0, 1]$ and $m>0$. In this connection, $\theta = 0$ case is already considered by D'Abbicco-Ikehata \cite{DI}, so one should restrict the parameter $\theta$ to the case $\theta > 0$. In this sense, this study is a kind of generalization of \cite{DI} to the general $\theta \in (0,1)$ through the logarithmic Laplacian type of mass-like term. Asymptotic behavior of the solution and its optimal decay/growth property has been already discussed in \cite{CAI, CAI-DCDS, JDE-2022} in the case when the logarithmic Laplacian with various physical sense can be included in the equation itself, and it motivates us to study (1.1) as one more topic. \\

We first define several function spaces related to the logarithmic Laplace operator $L_{\theta}$ introduced by Char\~ao-Ikehata \cite{Log-damping} (for $\theta=1$) and Char\~ao-D'Abbicco-Ikehata \cite{CDI} such that for $\delta \geq 0$
\begin{align}\label{Ydeltaspace}
Y_{\theta}^{\delta}  := \{ f \in L^2({\bf R}^n) :  \int _{{\bf R}_{\xi}^{n}} (1+\log (1+|\xi|^{2\theta} ))^{\delta} |\hat{f} (\xi)|^2d \xi < \infty  \}
\end{align}
with its natural norm
\begin{equation}\label{definormYdelta}
\Vert f\Vert_{Y_{\theta}^{\delta}} := \left(\int _{{\bf R}_{\xi}^{n}} (1+\log (1+|\xi|^{2\theta} ))^{\delta} |\hat{f} (\xi)|^2d \xi\right)^{1/2}, \quad f \in Y_{\theta}^\delta ,
\end{equation}
and its corresponding inner product.

\begin{rem}{\rm Due to the fact that  $\log (1+|\xi |^{2\theta}) \leq |\xi|$ for all $\xi \in {\bf R}^n$ satisfying $\vert\xi\vert \gg 1$,  one notices $ H^{s}({\bf R}^{n}) \subset Y_{\theta}^{s} \subset L^2$ for $s \geq 0$.}
\end{rem}

The linear operator  
\[L_{\theta}: D(L_{\theta}) \subset L^{2}({\bf R}^{n}) \to L^{2}({\bf R}^{n}), \quad \theta>0,\]
is defined as follows: 
\[D(L_{\theta}) := \left\{f \in L^{2}({\bf R}^{n}) \,\bigm|\,\int_{{\bf R}^{n}}(\log(1+\vert\xi\vert^{2\theta}))^{2}\vert\hat{f}(\xi)\vert^{2}d\xi < +\infty\right\} = Y_{\theta}^{2},\]
and for $f \in D(L_{\theta})$,  
\[(L_{\theta}f) (x) := {\cal F}_{\xi\to x}^{-1}\left(\log (1+\vert\xi\vert^{2\theta})\hat{f}(\xi)\right)(x).\]
\noindent
Here, one has just denoted the Fourier transform ${\cal F}_{x\to\xi}(f)(\xi)$ of $f(x)$ by 
\[{\cal F}_{x\to\xi}(f)(\xi) = \hat{f}(\xi) := \displaystyle{\int_{{\bf R}^{n}}}e^{-ix\cdot\xi}f(x)dx, \quad \xi \in {\bf R}^n,\]
as usual with $i := \sqrt{-1}$, and ${\cal F}_{\xi\to x}^{-1}$ expresses its inverse Fourier transform. 
\noindent

Since $L_{\theta}$ is constructed by a nonnegative-valued multiplication operator it is nonnegative and self-adjoint in $L^{2}({\bf R}^{n})$ (see \cite{Log-damping}), the square root 
$$L_{\theta}^{1/2}: D(L_{\theta}^{1/2}) \subset L^{2}({\bf R}^{n}) \to L^{2}({\bf R}^{n})$$
can be defined, and is also nonnegative and self-adjoint with its domain
\[D(L_{\theta}^{1/2}) = \left\{f \in L^{2}({\bf R}^{n}) \,\bigm|\,\int_{{\bf R}^{n}}\log(1+\vert\xi\vert^{2\theta })\vert\hat{f}(\xi)\vert^{2}d\xi < +\infty\right\} = Y_{\theta}^{1}.\] 
Note that $Y_{\theta}^{1}$ becomes Hilbert space with its graph norm
$$\Vert v\Vert_{D(L_{\theta}^{1/2})} := \left(\Vert v\Vert^{2} + \Vert L_{\theta}^{1/2}v\Vert^{2}\right)^{1/2},$$
where to simplify the notation we define the  $L^{2}({\bf R}^{n})$-norm  by
\[\Vert\cdot\Vert := \Vert\cdot\Vert_{L^{2}({\bf R}^{n})}.\]
\noindent
We  also note that  
$$H^{s}({\bf R}^{n}) \hookrightarrow Y_{\theta}^{1} \hookrightarrow  L^{2}({\bf R}^{n})$$
for any $ s > 0$. Symbolically writing, one can see
\[L_{\theta}= \log(I+(-\Delta)^{\theta}),\]
where $\Delta$ is the usual Laplace operator defined on $H^2({\bf R}^n)$. For more detailed and recent studies concerning the logarithmic Laplacian, one can cite \cite{B, B-2}, \cite{CW}, \cite{N} and \cite{WJ}.
 
The associated energy identity is given by 
\begin{equation}\label{energy}
E_{u}(t)  +  \int_0^t\Vert \nabla u_s(s,\cdot)\Vert^{2}ds=  E_{u}(0), \quad t>0,
\end{equation} 
where
\[
E_{u}(t) := \frac{1}{2}\left(\Vert u_{t}(t,\cdot)\Vert_{L^{2}}^{2} + \Vert  \nabla u(t,\cdot)\Vert_{L^{2}}^{2} + m^{2}\Vert L_{\theta}^{1/2}u(t,\cdot)\Vert_{L^{2}}^{2}\right).
\]
The identity \eqref{energy} implies that the total energy is a non-increasing function in time because of the existence of the  strong damping  term $  -\Delta u_{t}$ in \eqref{equation}.  \\

Let us recall several known historical facts related to our dissipative model (1.1), which has double dispersion.\\
First of all,  one should cite two celebrated works due to \cite{P} and \cite{S}, which studies the diffusive aspect and $L^{p}$-$L^{q}$ estimates of the solution to the strongly damped wave equation: 
\begin{equation}\label{Ryo-1}
u_{tt} -\Delta u -\Delta u_{t} = 0.
\end{equation}
After these works, so many precise analysis can be applied to \eqref{Ryo-1} through the papers due to \cite{I-14, IO, IT, ITY}, where the authors investigate the asymptotic profile and optimal decay and/or non-decay property of the solution itself. Precise decay estimates and its application to nonlinear problem can be discussed in \cite{DR}. $L^{1}$-estimates of the solution can be derived deeply in \cite{NR}. Higer order asymptotic expansions of the solution itself and the square of the $L^{2}$-norm of the solution as $t \to \infty$ can be well-investigated in \cite{BV-1, BV-2} and \cite{Mi}. Abstract theories concrning the regularity and decay property of \eqref{Ryo-1} can be constructed in \cite{GGH}, \cite{ITY} and \cite{LIC}.

Next, a generalization of \eqref{Ryo-1} to the so-called Klein-Gordon equation with strong damping such that
\begin{equation}\label{Ryo-2}
u_{tt} -\Delta u + m^{2}u -\Delta u_{t} = 0
\end{equation}
can be considered in \cite{DI}, and the effect of oscillation in the low frequency part of the solution captured to \eqref{Ryo-1} has been disappeared in the equation \eqref{Ryo-2}. This implies that the additional mass-like term $m^{2}u$ makes the structure of the solution itself change. Thus, it is quite natural to consider the following generalization of \eqref{Ryo-2}: 
\begin{equation}\label{Ryo-3}
u_{tt} + (-\Delta)^{\sigma}u + m^{2}(-\Delta)^{\theta}u + (-\Delta)^{\delta}u_{t} = 0,
\end{equation}
where $\sigma > 0$, $\theta > 0$ and $\delta > 0$. To the best of authors' knowledge, this type of equation is studied at least by two papers \cite{D} nd \cite{CD}. In fact, in the case of $m > 0$, $\sigma > \theta > 0$, and $\delta > 0$ the author in \cite{D} studies several decay rates of the solution and the anomalous diffusion phenomenon, while the authors in \cite{CD} treat the case of $m = 1$, $\sigma = 2$, $\theta = 1$, and $\delta =1$, and investigates $L^{p}$-$L^{q}$ estimates fo $1 \leq p \leq q \leq \infty$ of the solution, and inviscid limits for the first and the second order profiles together with some nonlinear effect, however, it seems that these authors are not conscious to observe the growth property of the solution itself at this stage. 

In this connection, in the case when $m > 0$, $\sigma = 1$, $\theta = 0$, and $\delta = 0$ in \cite{Z} the author studied the exponential decay of the corresponding total energy $E_{m}(t)$ such that
\[E_{m}(t) := \frac{1}{2}\left(\Vert u_{t}(t,\cdot)\Vert^{2} + \Vert \nabla u(t,\cdot)\Vert^{2} + m^{2}\Vert u(t,\cdot)\Vert^{2}  \right).\]
In fact, the author of \cite{Z} considered the nonlinear problem in a more general framework.\\ 

Under these previous developments it seems to be little investigations on the topic to observe the growth property (as $t \to \infty$) of the solution itself. Our main concern is to consider the equation (1.1) by replacing the mass-like term $m^{2}(-\Delta)^{\theta}u$ of \eqref{Ryo-3} with $\sigma = \delta = 1$ by logarithmic mass-like one $m^{2}\log(I + (-\Delta)^{\theta})u$ to soften the regularity of the solution, and to capture the growth property of the solution itself coming from an oscillation effect in the low frequency region of the solution. In this situation, one wants to observe the effect of the mass-like term brought by a new one $\log(I + (-\Delta)^{\theta})$ on the decay structure of the equation. To catch information on the time behavior of the $L^{2}$-nom of the solution itself is very important to observe some singularity (near low frequency part) included in the solution itself to problem (1.1)-(1.2). Such a singularity reflects on the quantity of the function $t \mapsto \Vert u(t,\cdot)\Vert$, not on the time and/or spatial derivatives of the solutions such as total energy.\\  

Now, we introduce the asymptotic profile as $t \to \infty$ of the solution to problem (1.1)-(1.2):
\begin{equation}\label{Perfil-1}
\varphi(t,\xi)=P_1e^{-\frac{\vert\xi\vert^{2}}{2}t} \frac{\sin\big(t\sqrt{|\xi|^2 +m^2\log(1+|\xi|^{2\theta} )} \big)}{\sqrt{|\xi|^2+m^2 \log(1+|\xi|^{2\theta})}}, 
\end{equation}
where the moment of the initial velocity $P_{1} \in {\bf R}$ is defined by 
\[P_{1} := \int_{{\bf R}^{n}}u_{1}(x)dx.\]
Then, our main result reads as follows.
\begin{theo}\label{Theorem1}
Let $n \geq 1$, and $0 < \theta \leq 1$. Let $\varphi$ be the function defined in \eqref{Perfil-1}. Choose $u_{0} = 0$, and $u_{1} \in L^{1,\theta}({\bf R}^{n})\cap L^{2}({\bf R}^{n})$. Then there exists a positive constant $C=C_{n,\delta_0,\theta}$ such that the mild solution $u \in C([0,\infty);H^{1}({\bf R}^{n})) \cap C^{1}([0,\infty);L^{2}({\bf R}^{n}))$ to problem {\rm (1.1)-(1.2)} satisfies
\begin{align*}
\int_{{\bf R}^n} |\hat{u}(t,\xi) - \varphi(t ,\xi)|^2d\xi 
%& \leq C\Big( \sum_{i=1}^3 \int_{|\xi|\leq \delta_0}|F_i(t,\xi)^2d\xi 
%+  \int_{|\xi|\geq \delta_0}|\hat{u}(t,\xi)^2d\xi  +  \int_{|\xi|\geq \delta_0}|\varphi(t,\xi)^2d\xi \Big)\\
&\leq C_{n,\delta_0,\theta}\Big( m^{-2}P^2_1 t^{-\frac{n+8-6\theta}{2}}
+m^{-2}P^2_1 t^{-\frac{n+4-4\theta}{2}} + m^{-2}||u_1||^2_{1,\theta} t^{-\frac{n}{2}} + P_{1}^{2}e^{-\gamma t}\\
&+ ||u_1||^2_{1} t^{2}e^{-\gamma t} + ||u_1||^2_{1} t^{2}e^{-\alpha t}
+ ||u_1||^2_{1} t^{2}e^{-\beta t}
\Big), \quad t \geq 0, \quad 
%%%% 
%%%
\end{align*} 
where $\alpha > 0$, $\beta > 0$ and $\gamma > 0$ are generous constants. 
\end{theo}

As a reult of Theorem \ref{Theorem1}, one can state the following optimal $L^{2}$-decay and/or growth estimates of the solution itself to problem (1.1)-(1.2) 
\begin{theo}\label{Theo6.2}\,Under the same assumptions as in Theorem {\rm 1.1} the following statements are true: 
\begin{itemize}
\item[{\rm (i)}] for $ n=1$ with $ 0< \theta < \frac{1}{2} $, $n=2$ with $ 0< \theta < 1 $ and $n \geq 3  $ with $ 0<  \theta \leq 1 $, there exist positive constants $C_1, C_2$ depending on $n$, $\theta$ and $m > 0$ such that
$$C_1 |P_1| t^{-\frac{n-2\theta}{4}}\leq \|u(t, \cdot )\| \leq C_2  \Big ( \Vert u_{1}\Vert_{1} + \|u_1\|_{1,\theta} \Big ) t^{-\frac{n-2\theta}{4}}, \quad t\gg 1,$$
\item[{\rm (ii)}] for $n=1 $ and $ \frac{1}{2}< \theta \leq 1$, there exist positive constants $C_1, C_2$ such that  
$$ \frac{C_1}{(1+m^{2})^{\frac{1}{4\theta}}} |P_1| t^{ \frac{2\theta -1}{2\theta} } \leq \| u(t, \cdot ) \| \leq \frac{C_2}{m} \frac{1}{\sqrt{2\theta -1}} \Big ( \Vert u_{1}\Vert_{1}+ \|u_1\|_{1,\theta} \Big )  t^{ \frac{2\theta -1}{2\theta} } , \quad t \gg 1,$$  
\item[{\rm (iii)}] for $n=1$ and $\theta =\frac{1}{2}$ or $n = 2$ and $\theta = 1$, there exist positive constants $C_1, C_2$ such that  
$$\frac{C_1}{\sqrt{2+m^{2}}} |P_1| \sqrt{\log t} \leq \| u(t, \cdot ) \| \leq \frac{C_2}{m} \Big ( \Vert u_{1}\Vert_{1} + \|u_1\|_{1,\theta} \Big ) \sqrt{\log t} , \quad t \gg 1.$$  
%$$n=2, \theta = 1 \Rightarrow \frac{C_1}{\sqrt{2+m^{2}}} |P_1| \sqrt{\log t} \leq \| u(t, \cdot ) \| \leq \frac{C_2}{m} \Big ( |P_1|+ \|u_1\|_{1,\theta} \Big ) \sqrt{\log t} , \quad t \gg 1. $$  
\end{itemize} 
\end{theo}

Let us reconsider the role of mass-like term from Theorem 1.2.
\noindent
The result of (i) in Theorem 1.2 implies the decay estimate, and this comes from the stronger effect with decay order $t^{-\frac{n}{4}}$ of the Gauss kernel $e^{-t\vert\xi\vert^{2}/2}$ in the low frequency zone, and the oscillation effect comimg from the mass-like term $m^{2}\log(1+(-\Delta)^{\theta})u \sim m^{2}(-\Delta)^{\theta}u$ gives a small role with growth order $t^{\frac{\theta}{2}}$. Totally, one can get the optimal decay rate $t^{-\frac{n}{4} + \frac{\theta}{2}} = t^{-\frac{n-2\theta}{4}}$. Recall that in case of $n \geq 3$ the optimal decay order of the $L^{2}$-norm of the solution to \eqref{Ryo-1} is $t^{-\frac{n-2}{2}}$ (see \cite{I-14}). In this case, one can see that $\frac{n-2}{2} < \frac{n-2\theta}{2}$ for $\theta \in [0,1]$. In particular, in the case of $\theta = 0$, the decay order coincides with that derived in \cite{DI}, and the result corresponds to the strongly damped Klein-Gordon equation \eqref{Ryo-2}. Thus, one can understand that the additional term $m^{2}\log(1+(-\Delta)^{\theta})u$ plays a role to make the decay rate more fast. In this sense, the name of "mass-like" term seems suitable. 
While, (ii) of Theorem 1.2 reflects a stronger effect with growth order $t$ of the oscillation factor in the low frequency part: 
\[\frac{\sin\big(t\sqrt{|\xi|^2 +m^2\log(1+|\xi|^{2\theta} )} \big)}{\sqrt{|\xi|^2+m^2 \log(1+|\xi|^{2\theta})}} \sim \frac{\sin(mtr^{\theta})}{mr^{\theta}} \quad (r \to +0),\]
and the diffusion part is no-effective. Roughly observing, for $n = 1$ one can see that
\[\left(\int_{0}^{1}\frac{\sin^{2}(mtr^{\theta})}{(mr^{\theta})^{2}}dr\right)^{1/2} \sim t\left(\int_{0}^{m^{-\frac{1}{\theta}}t^{-\frac{1}{\theta}}}dr\right)^{1/2} \sim t\cdot t^{-\frac{1}{2\theta}} \quad (t \to \infty),\]
since $\vert\displaystyle{\frac{\sin\eta}{\eta}}\vert \sim 1$ when $\eta \to +0$. As a result one can get infinite time blow up rate $t\cdot t^{-\frac{1}{2\theta}} = t^{\frac{2\theta-1}{2\theta}}$ when $t \to \infty$ as shown in (ii). In some sense, in (ii) of Theorem 1.2 the dominant part in the low frequency zone as $t \to \infty$ of the equation (1.1) will be
\[u_{tt} + m^2(-\Delta)^{\theta}u = 0.\]
(iii) of Thorem 1.2 is just an intermediate case between (i) and (ii). These reflect a removable singularity in the low frequency zone included in the solution itself. This kind of growth property as shown in (ii) and (iii) has been newly discovered to the equation (1.1) itself.
In this connection, a regularity aspect which corresponds to high frequency asymptotics of the solution to (1.1) can be dominated by \eqref{Ryo-1} even if $\theta > 0$ is large enough (cf., \cite{CD}).\\

{\bf Notation.} {\small Throughout this paper, $\| \cdot\|_q$ stands for the usual $L^q({\bf R}^{n})$-norm. For simplicity of notation, in particular, we use $\| \cdot\|$ instead of $\| \cdot\|_2$. By $(\cdot,\cdot)$ we mean the usual $L^{2}$-inner product. Furthermore, we denote $\Vert\cdot\Vert_{H^{l}}$ as the usual $H^{l}({\bf R}^{n})$-norm. We  also define a relation $f(t) \sim g(t)$ as $t \to \infty$ by: there exist constant $C_{j} > 0$ ($j = 1,2$) such that
\[C_{1}g(t) \leq f(t) \leq C_{2}g(t)\quad (t \gg 1).\] 

For $\Omega \subset {\bf R}^n$ we denote $f \approx g$  on $\Omega$, if and only if 
there are constants $K_1, K_2$ such that 
$$ K_1 f(y) \leq g(y) \leq K_2 f(y), \;\;\mbox{for all} \; y \in \Omega.$$

We also introduce the following weighted functional spaces for $\gamma >  0$:
\[L^{1,\gamma}({\bf R}^{n}) := \left\{f \in L^{1}({\bf R}^{n}) \; \bigm| \; \Vert f\Vert_{L^{1,\gamma}} := \int_{{\bf R}^{n}}(1+\vert x\vert^{\gamma})\vert f(x)\vert dx < +\infty\right\}.\]
${\rm Re}z$ and ${\rm Im}z$ denote the real and the imaginary parts of $z$, respectively. Finally, we denote the surface area of the $n$-dimensional unit ball by $\omega_{n} := \displaystyle{\int_{\vert\omega\vert = 1}}d\omega$. 

}

%%%%%%%%%%%%%%%%%%E½@PROOFE½@%%%%%%%%%%%%%%%%%%%%%%%%
%%%%%%%%%%%%%%%%%%E

\section{Basic preliminary results}
In this section we shall collect important lemmas to derive precise estimates of the several  quantities related to the solution to problem (1.1)-\eqref{initial}. These are already studied and developed in our previous works (see \cite{Log-damping, CDI}).

%The following estimate for the function  
%$$I_p(t)= \int_0^{1}(1+r^{2})^{-t}r^p dr$$ 
%is a direct consequence of the cases $p \geq 0$ in Char\~ao-Ikehata \cite{Log-damping} and $-1<p<0$ in Char\~ao-D'Abbicco-Ikehata \cite{CDI}.

%\begin{lem}\label{general-p}
% Let $p > -1$ be a real number. Then
%$$I_p(t) \sim t^{-\frac{p+1}{2}}, \quad t \gg 1.$$
%%
%\end{lem}

%In order to deal with the high frequency part of estimates, one defines a function again 
%$$J_p(t)=\int_1^{\infty}(1+r^2)^{-t}r^p dr$$
%for $p \in {\bf R}$.

%Then the next lemma is important to get estimates on the zone of high frequency to the solutions of the problem \eqref{eqn}--\eqref{initial}. The proof  appears in Char\~ao-Ikehata \cite{Log-damping}.
%\begin{lem}\label{infit}
%\,Let $p \in {\bf R}$. Then it holds that 
%$$J_p(t) \sim \dfrac{2^{-t}}{t-1}, \quad t \gg 1.$$
%\end{lem}
%\vspace{0.2cm}
%For later use we prepare the following simple lemma, which implies the exponential decay estimates of the middle frequency part.
%\begin{lem}\label{intermid}\,Let $p \in {\bf R}$, and $\eta \in (0,1]$. Then there is a constant $C > 0$ such that 
%$$\int_{\eta}^{1}(1+r^{2})^{-t}r^{p}dr \leq C(1+\eta^{2})^{-t}, \quad t \geq 0.$$
%\end{lem}

%%%%%%%%%%%%%%%%%%%%%%%%%%%%%%%%%%%%%%%%%%%%%
%\begin{rem}
%{\rm We note that the proof of Lemma \ref{general-p} is done by using simple differential calculus and the theory from hypergeometric functions (see Watson \cite{W}). These are already developed in \cite{Log-damping} and \cite{CDI}.}
%\end{rem}
%%%%%%%%%%%%%%%%%%%%%%%%%%%%%%%%%%%%%%

The following two lemmas are well known.
\begin{lem}\label{lem-t}
Let  $k>-n$, $\alpha> 0$ and  $C>0$. Then there exists a constant  $K>0$ depending only  on  $n, \alpha, k$ such that
\begin{align*}
\int_{|\xi|\leq 1} e^{-C|\xi|^{\alpha}t}|\xi|^k d\xi = \omega_{n}\int_{0}^{1}e^{-Cr^{\alpha}t}r^{n-1+k}dr \leq K t^{-\frac{n+k}{\alpha}},
\end{align*}
for all $t \geq 1$.
\end{lem} 

\begin{lem}\label{lem-sen-hiper}
$$ \frac{\sinh x }{x} \leq   e^{x}  $$
for all $ x >0 $. 
\end{lem}

One also needs the following decomposition of the Fourier transform of a function $f \in L^1({\bf R}^n)$  as follows:   
\begin{equation} \label{decompo}
\hat{f}(\xi)=A_f(\xi)-iB_f(\xi)+P_f,
\end{equation}
for all $\xi \in {\bf R}^n$ , where 
\begin{itemize}
\item[$\bullet$] $A_f(\xi)=\displaystyle{\int_{{\bf R}^n}}{(\cos(x\cdot\xi)-1)f(x)}dx,$
\item[$\bullet$] $B_f(\xi)=\displaystyle{\int_{{\bf R}^n}}{\sin(x\cdot\xi)f(x)}dx,$
\item[$\bullet$] $P_f=\displaystyle{\int_{{\bf R}^n}}{f(x)}dx.$
\end{itemize}
\noindent
Then, the next lemma has been prepared in \cite{I-04} (see Notation for the definition of $L^{1,\kappa}({\bf R}^{n})$).
\begin{lem}\label{lema2.6}
\begin{itemize}
\item[{\rm i)}] If\;  $f \in L^1({\bf R}^n)$  then for all $\xi \in {\bf R}^n$ it is true that  
$$|A_f(\xi)|\leq L\|f\|_{L^1} \quad \text{ and  } \quad |B_f(\xi)|\leq N\|f\|_{L^1}.$$
\item[\rm{ii)}] If \;$0<\kappa \leq 1$ and  $f \in  L^{1,\kappa}({\bf R}^n)$ then for all  $\xi \in {\bf R}^n$ it is true that 
$$|A_f(\xi)|\leq K|\xi|^\kappa\|f\|_{L^{1,\kappa}} \quad \text{ and  } \quad |B_f(\xi)|\leq M|\xi|^\kappa\|f\|_{L^{1,\kappa}}$$
\end{itemize}
\noindent with  $L$, $N$, $K$ and  $M$ positive  constants  depending only on  the dimension  $n$  or $n$ and $\kappa$.\\
\end{lem}

In order to justify the existence of a unique solution to problem
\eqref{equation}-\eqref{initial} we prepare the following theorem (cf. \cite[Theorem 6.4]{Goldstein}).

%\begin{theo}\label{Theo6.1Goldstein}{\rm [\cite[6.1, Goldstein]{Goldstein}]}
%Let $A$ be a generator of  contraction $C^0$-semigroup in $X$ and $B$ a dissipative operator with domain $D(B)$ such that $ D(A) \subset D(B) $ in $X$. Assume that there exist constants $ a, b \geq 0 $ with $ a <1  $ such that  
%$$ \| Bf \|_X \leq a \| Af \|_X + b \| f \|_X , \quad \forall  f \in D(A) .$$ 
%Then, $ A+B : D(A) \rightarrow X $ is a generator of a $C^0$-semigroup of contraction. 
%\end{theo}

%\begin{theo}\label{Theo6.4,Goldstein}{\rm [\cite[Theorem 4.6]{Goldstein}]} Let $A$ be a generator of $C^0$-semigroup in $X$ and $B \in \mathcal{B}(X) $, that is, $B: X \rightarrow X$ is a bounded linear operator. Then $ A+B$ generates a $C^0$-semigroup. 
%\end{theo}

\begin{theo}\label{AdaptTheorem} Let  $X$ be a Hilbert space and $A \subset X: D(A) \rightarrow X$ be a generator of a $C^0$-semigroup in $X$. If a linear operator $B: X \rightarrow X $ is bounded, then the operator $A+B$ generates a $C^0$-semigroup on $X$. 
\end{theo}
%{\it{Proof.}}
%The proof is a simple adaptation of the one presented for the Theorem 6.4 in \cite{Goldstein}.  In fact, we consider the operator $ C:= B- \|B\|I : D(B) \rightarrow X $ (here the symbol $\| B \|$ means the norm of the operator $B$ and $I$ is the identity operator). We can see the operator $C$ is dissipative and 
%\begin{align*}
%\|C f \|_X = \| (B-\|B\|I) f \|_X \leq \|B f\| + \|B\| \|f\| \leq 2\|B\| \|f\|.
%\end{align*}

%Then  we may conclude, from Theorem \ref{Theo6.1Goldstein}, that the operator $A+C$ generates a $C^0$-semigroup of contractions. Since $ \|B\|I \in \mathcal{B}(X) $, the operator $$ A+B = A+C +\|B\|I   $$
%generates a $C^0$-semigroup according to Theorem \ref{Theo6.4,Goldstein}. 
%\hfill
%$\Box$

%%%%%%% We don't use the following lemma %%%%%
%The proof of the following lemma is now standard.
%\begin{lem}\label{lemma2.10} Let $n \geq 1$, and $k > -n $, $\nu>0$ and $C>0 $. Then there exists  a %constant $K>0$ which depends only on $n$ such that 
%$$ \int _{{\bf R}^n} e^{-C|\xi|^{\nu}t}|\xi|^{k}d\xi \leq Kt^{- \frac{n+k}{\nu}}$$
%for all $t>0$.
%\end{lem}
%{\it Proof.} The proof ot this lemma appears in \cite{....................}.
%\hfill
%$\Box$

%%%%%%%%%%%%%%%%%%%%%%%%%%%%%%%%%%%%%%%%%%%%%%%%%%%%%%%%%%%%%%%%%%%%%%%%%%%%%%%%%%%%%%%%%%%%%%%%%%%%%%%%%%%%%%%%%%%%%%%%%%%%%%%%%%%%%%%%

\section{Existence and Uniqueness }

 The total energy associated to the equation \eqref{equation} is defined  by
\begin{align*}
E(t):= \frac{\| u_t \|^2 + \| \nabla u \|^2 + \|L_\theta ^{\frac{1}{2}} u \|^2}{2}, \quad t \geq 0, 
\end{align*} 
and it satisfies the following energy identity 
\begin{align*}
\frac{d}{dt} E(t) + \| \nabla u_t \|^2 =0, \quad t>0.
\end{align*}
So it is natural to consider the following energy space
\begin{align}
X= (H^1({\bf R}^{n}) \cap D(L_\theta^{1/2})) \times L^2. 
\end{align}
Here, one can observe that $ H^1({\bf R}^{n}) \subset D(L_\theta^{1/2}))  $ since $ 0\leq \theta \leq 1 $, and due to the facts that $0\leq \log(1+r^{2\theta}) \leq \log 2 <1 $ for $ 0\leq r \leq 1$, and because of $ \log(1+r^{2\theta}) \leq r^2 $  for $r\geq 1$.  Thus, in the case $ 0\leq \theta \leq 1 $ the energy space can be sufficiently chosen as 
\begin{align}
X= H^1({\bf R}^{n}) \times L^2({\bf R}^{n}). 
\end{align}
We remember that the space $ H^1({\bf R}^{n})  $ is a Hilbert space with the natural norm given by 
$$ \| f \| _{H^1} :=  \left( \int _{{\bf R}^{n}} (1+ |\xi|)^2 |\hat{f}(\xi)|^2 d\xi \right) ^{\frac{1}{2}},$$ 
but we may observe that this norm is equivalent to the following norm
$$ \|| f \|| _{{H^1}} :=  \left( \int _{{\bf R}^{n}} (1+ m^2 \log (1+|\xi|^{2\theta}) +|\xi|^2) |\hat{f}(\xi)|^2 d\xi \right) ^{\frac{1}{2}} , \quad f \in H^1({\bf R}^n) ,$$
since $0\leq \theta \leq 1 $. This second norm can be defined by the following inner product
\begin{equation}\label{P4innerproduct}
 ((f,g))_{H^1}:= \int _{{\bf R}^{n}} (1+ m^2 \log (1+|\xi|^{2\theta}) +|\xi|^2) \hat{f}(\xi) \overline{\hat{g}(\xi)}d\xi  . 
\end{equation} 
Therefore, in this section we consider the Hilbert space $ H^1({\bf R}^{n}) $ with the inner product $  ((\cdot,\cdot))_{H^1} $ just defined in \eqref{P4innerproduct}. 
\vspace{0.2cm}
In order to prove the existence and uniqueness, we reduce the original equation \eqref{equation} to the first order system as usual.\\
We set $ v= u_t $, and $ U = \begin{pmatrix}
u \\ v
\end{pmatrix}  $ and, we have (at least formally) 
\begin{align}
\frac{\mathrm{d} }{\mathrm{d} t} U &=  \frac{\mathrm{d} }{\mathrm{d} t} \begin{pmatrix}
u \\ v
\end{pmatrix} = \frac{\mathrm{d} }{\mathrm{d} t} \begin{pmatrix}
u \\ u_{t}
\end{pmatrix}  =  \begin{pmatrix}
u_t \\
\Delta u_t + \Delta u -m^2 L_\theta u 
\end{pmatrix}  \nonumber \\
&=  \begin{pmatrix}
v \\
\Delta (v + u) -m^2 L_\theta u 
\end{pmatrix}\nonumber\\
%= \begin{pmatrix}
%0 & I \\
%\Delta - m^2 L_\theta  & \Delta  \\
%\end{pmatrix}
%\begin{pmatrix}
%u \\ v
%\end{pmatrix} \nonumber \\
&= \begin{pmatrix}
0 & I \\
\Delta - m^2 L_\theta - I  & \Delta  \\
\end{pmatrix} \begin{pmatrix}
u \\ v
\end{pmatrix} + \begin{pmatrix}
0 \\ u
\end{pmatrix} \nonumber \\
&=: BU + FU,
\end{align}

%In order to prove the existence and uniqueness, we reduce the original equation \eqref{equation} to the first order system as usual.\\
%We set $ v= u_t $, and $ U = \begin{pmatrix}
%u \\ v
%\end{pmatrix}  $ and, we have (at least formally) 
%\begin{align}
%\frac{\mathrm{d} }{\mathrm{d} t} U &=  \frac{\mathrm{d} }{\mathrm{d} t} \begin{pmatrix}
%u \\ v
%\end{pmatrix} = \frac{\mathrm{d} }{\mathrm{d} t} \begin{pmatrix}
%u \\ u_{t}
%\end{pmatrix}  =  \begin{pmatrix}
%u_t \\
%\Delta u_t + \Delta u -m^2 L_\theta u 
%\end{pmatrix}  \nonumber \\
%&=  \begin{pmatrix}
%v \\
%\Delta v + \Delta u -m^2 L_\theta u 
%\end{pmatrix}  = \begin{pmatrix}
%0 & I \\
%\Delta - m^2 L_\theta  & \Delta  \\
%\end{pmatrix} \begin{pmatrix}
%u \\ v
%\end{pmatrix} \nonumber \\
%&= \begin{pmatrix}
%-(-\Delta)^{1/2} & I \\
%\Delta - m^2 L_\theta - I  & \Delta  \\
%\end{pmatrix} \begin{pmatrix}
%u \\ v
%\end{pmatrix} + \begin{pmatrix}
%(-\Delta)^{1/2}u \\ u
%\end{pmatrix} \nonumber \\
%&=: BU + FU,
%\end{align}
%where
%$B$ and $F$ are given by 

%\begin{align}\label{P4.DefB}
%BU = \begin{pmatrix}
%-(-\Delta)^{1/2} & I \\
%\Delta - m^2 L_\theta - I  & \Delta  \\
%\end{pmatrix} U 
%\end{align}
where
\begin{align}\label{P4.DefB}
BU := \begin{pmatrix}
v\\ \Delta(u+v) - m^{2}L_{\theta}u - u
\end{pmatrix} .
\end{align}
for $U := (u,v) \in D(B)$,
and 
\begin{align}\label{P4.DefF}
FU = \begin{pmatrix}
0 \\ u
\end{pmatrix}
\end{align}
for  $U := (u,v) \in X$.
\noindent
In this case, we consider the domain $D(B)$  of the linear operator $B$ as  
\begin{align}\label{B-domain}
D(B)= \{ (u,v) \in  H^1({\bf R}^n) \times H^1({\bf R}^n) \;\;   \mbox{such that } u+v \in H^2({\bf R}^n) \}.
\end{align}
This idea is coming from \cite{ITY}, and it is crucial to prove Lemma 3.2. The domain $D(B)$ is dense in the energy space $X=H^1({\bf R}^n) \times L^2({\bf R}^n)$. 

Let us see that the linear operator $B$ is dissipative in $X$, and satisfies $ (I-B)(D(B))=X $.

\begin{lem} \label{B-dissip}
The linear operator $B$ is dissipative in $X$. 
\end{lem}
{\it{Proof.}} In the space $X=H^1({\bf R}^n) \times L^2({\bf R}^n)$, we consider the natural inner product
$$ \left ( (u_1,v_1),(u_2, v_2) \right )_{X}= ((u_1,u_2))_{H^1}+ (v_1,v_2) $$ 
for $(u_j, v_j) \in H^1({\bf R}^n) \times L^2({\bf R}^n)$, ($j=1,2 $), where  $ (( \cdot , \cdot))_{H^1} $ is given by \eqref{P4innerproduct} and $(\cdot, \cdot)$ denotes the usual inner product in $L^2$. 

Let $(u,v) \in D(B)$. Then, the definition of the equivalent inner product \eqref{P4innerproduct} in $H^1({\bf R}^n)$ implies that

%\[\left(( B(u,v),(u,v))\right )_X = \left(( -(-\Delta)^{1/2}u + v, u)\right)_{H^{1}} + (\Delta(u+v)-m^{2}L_{\theta}u - u, v)\]
%\[= \int_{{\bf R}^{n}}A(r)(-\vert\xi\vert\hat{u} + \hat{v})\bar{\hat{u}}d\xi - \int_{{\bf R}^{n}}\vert\xi\vert(\hat{u}+\hat{v})(\vert\xi\vert\bar{\hat{v}})d\xi - \int_{{\bf %R}^{n}}\left(m^{2}\log(1+\vert\xi\vert^{2\theta})\hat{u} + \hat{u}\right)\bar{\hat{v}}d\xi\]%
%\[= 2i\int_{{\bf R}^{n}}A(r){\rm Im}(\hat{v}\bar{\hat{u}})d\xi - \int_{{\bf R}^{n}}A(r)\vert\xi\vert\vert\hat{u}\vert^{2}d\xi - \int_{{\bf R}^{n}}\vert\xi\vert^{2}\vert\hat{v}\vert^{2}d\xi,\]
%where
%\[A(r) := 1+ m^{2}\log(1+\vert\xi\vert^{2\theta}) + \vert\xi\vert^{2} = 1+ m^{2}\log(1+r^{2\theta}) + r^{2},\quad (r = \vert\xi\vert).\]
%\noindent
\begin{align*}
\left ( B(u,v),(u,v)  \right )_X&= \left((v, \Delta (u+v) -m^2L_{\theta}u-u), (u,v) \right)_{X} \\
&= ((v,u))_{H^1} + (\Delta (u+v) -m^2L_{\theta}u-u, v)
 \\
&= \int _{{\bf R}^n} (1+m^2\log(1+|\xi|^{2\theta})+|\xi|^2) \hat{v} \overline{\hat{u}} d\xi \\
&- \int _{{\bf R}^n} (1+m^2\log(1+|\xi|^{2\theta})+|\xi|^2) \hat{u} \overline{\hat{v}} d\xi - \int _{{\bf R}^n} |\xi|^2 |\hat{v}|^2 d\xi\\
&=2i\int _{{\bf R}^n} (1+m^2\log(1+|\xi|^{2\theta})+|\xi|^2) \text{Im}(\hat{v} \overline{\hat{u}} )d\xi - \int_{{\bf R}^n} |\xi|^2 |\hat{v}|^2 d\xi.
\end{align*}
Thus, one can get
$$ \text{Re}\left ( B(u,v),(u,v)  \right )_X =  - \int_{{\bf R}^{n}} |\xi|^2 |\hat{v}|^2 d\xi \leq 0, $$
and this concludes the desired statement. 
\hfill
$\Box$

\begin{lem}\label{P4.Lemmasurjective}
$ (I-B)(D(B))=X $. 
\end{lem}
 %Let $(u,v) \in D(B)$. Then,  the definition of $D(B)$ in \eqref{B-domain}, says that 
%\begin{align*}
%&u-v \in H^1(\R^n) \\
%&\Delta u-m^2\log(I+(-\Delta)^{2\theta})u - u+\Delta v = \Delta(u-v)
%-m^2\log(I+(-\Delta)^{2\theta})u - u  \in L^2(\R^n),
%\end{align*}
%because $u+v \in H^2(\R^n)$.
{\it{Proof.}}\,At first, the inclusion $(I-B)D(B) \subset X $ is easily checked. So, it suffices to check $X \subset (I-B)(D(B))$. 
%\vspace{0.2cm}
For this, let $ (f,g) \in X=H^1({\bf R}^{n}) \times L^2({\bf R}^{n}) $. We need to show that there exists a pair $(u,v) \in D(B) $ that satisfies the following two identities
%%5
\begin{align}
&\;u-v = f, \label{P4.Eq4.12} \\
&-\Delta(u + v) +  m^2L_{\theta}u +u +v  = g. \label{P4.Eq4.13}
\end{align}

Since we need to find $u,v$ such that 
$v = u-f$ and $u$ satisfying \eqref{P4.Eq4.13}
it is sufficient to obtain $u \in H^1$ such that 
\begin{align} \label{eq-B-sobre}
-2\Delta u +2u +  m^2L_{\theta}u=g-\Delta f +f.
\end{align}

Note that $g-\Delta f +f \in H^{-1}({\bf R}^{n})$ at this stage. Then we  may apply the Lax-Milgram Theorem to obtain the existence of a unique function $u \in H^1({\bf R}^{n})$ satisfying \eqref{eq-B-sobre} in weak sense. Since $v=u-f$ and $f\in H^1({\bf R}^{n})$ it follows that  $(u,v) \in H^1({\bf R}^{n}) \times H^1({\bf R}^{n})$.

Now from \eqref{P4.Eq4.13} we have 
$$-\Delta(u+v)  +u +v  = g- m^2L_{\theta}u \in L^2({\bf R}^{n})
$$
due to the assumption on $g$ and $u\in H^1({\bf R}^{n})$. Thus,  the elliptic regularity theorem says that $u+v \in H^2({\bf R}^{n})$.
Therefore $(u,v) \in D(B)$ and the operator $I-B$ is surjective. 
\hfill
$\Box$

\begin{theo}
$B: D(B) \rightarrow X$ is an infinitesimal generator of a $C^0$-semigroup of contractions. 
\end{theo}
{\it{Proof.}} It follows from the facts that  $D(B)$ is dense in $X$, $B$ is dissipative according to Lemma \ref{B-dissip},  and $R(I-B) = X$ by Lemma \ref{P4.Lemmasurjective}. The result follows from the Lumer-Phillips theorem.

\hfill
$\Box$

In order to conclude the existence and uniqueness of solutions to problem \eqref{equation}-\eqref{initial},  we observe that the operator $F$ given in  \eqref{P4.DefF} with domain $$ D(F) = H^1({\bf R}^{n}) \times L^2({\bf R}^{n}) $$ satisfies the following fact.
\begin{lem}
 $F: X \rightarrow X$ is a bounded linear operator. 
\end{lem}
{\it{Proof.}} Clearly the operator $F$ is linear. Let $ (u,v) \in X= H^{1}({\bf R}^{n}) \times L^2({\bf R}^{n}) $. Then, 
\begin{align*}
\| F(u,v)\|_X &= \|(0, u ) \|_{X} =   \|u\|  \leq \|(u,v)\|_{X},
\quad (u,v) \in X.
\end{align*}
Therefore, the operator $F$ is bounded in $X$. 
\hfill
$\Box$

Since $ B: D(B) \rightarrow X$ is an infinitesimal generator of a $C^0$-semigroup of contractions and $F: D(F) \rightarrow X$ is a bounded linear operator with $ D(B) \subset D(F)= X$, we may conclude the following result from Theorem \ref{AdaptTheorem}. 

\begin{theo} $B+F: D(B) (\subset X) \rightarrow X $ is an infinitesimal generator of a $C^0$-semigroup in $X$. 
\end{theo}

As a result, one can arrive at the following theorem concerning the existence and uniqueness to problem \eqref{equation}-\eqref{initial}. 

\begin{theo} Let $ n \geq 1 $  and  $ (u_0, u_1) \in  D(B)$, with $D(B)$ given by \eqref{B-domain}   Then the problem \eqref{equation}-\eqref{initial} admits a unique solution $u$ in the class 
$$ C^1([0,\infty), H^1({\bf R}^{n})) \cap  C^2([0,\infty), L^2({\bf R}^{n})) .$$
 Moreover, for initial data $ (u_0,u_1) \in  H^1({\bf R}^{n}) \times L^2({\bf R}^{n})  $ the problem admits a unique weak solution in the class 
$$ C([0,\infty), H^1({\bf R}^{n}))\cap  C^1([0,\infty), L^2({\bf R}^{n})).$$
\end{theo}

\section{Cauchy problem in the Fourier space }

The associated  Cauchy problem in the Fourier space to problem \eqref{equation}-\eqref{initial}, is given by 
\begin{align}\label{F-eq}
&\hat{u}_{tt} + |\xi|^2 \hat{u} +m^2 \log (1+|\xi|^{2\theta}) \hat{u} +|\xi|^2 \hat{u}_t = 0, \quad t >0, \quad \xi \in {\bf R}^n \\
& \hat{u}(0,\xi) = \hat{u}_{0}, \quad  \hat{u}_t(0,\xi)=\hat{u}_{1}, \quad  \xi \in {\bf R}^n. \label{F-ini}
\end{align}

The characteristic equation associated to problem \eqref{F-eq}--\eqref{F-ini} is given by 
$$ \lambda ^2 + |\xi|^2 \lambda + (|\xi|^2 +m^2 \log(1+|\xi|^{2\theta}))=0,$$ 
and its characteristics roots are
$$ \lambda _{\pm} = \frac{-|\xi|^2 \pm \sqrt{|\xi|^4 - 4 (|\xi|^2 + m^2 \log(1+|\xi|^{2\theta})) }}{2} .$$ 

We need to analyse the sign of the expression
$|\xi|^4 - 4 (|\xi|^2 + m^2 \log(1+|\xi|^{2\theta})) $, for $ \xi \in {\bf R}^n$
to decide whether the characteristics roots are real or complex.

To do that we consider the function
\begin{equation}\label{f-r}
 f(r) = r^4 - 4r^2  -4 m^2 \log(1+r^{2\theta}), \quad  r \geq 0,
\end{equation}
where $\theta \geq 0$. We set
\[g(r) := r^{2}(r^{2}-4),\quad h(r) := k\log(1+r^{2\theta}),\]
where $k := 4m^{2} > 0$. Then, one has
\[f(r) = g(r) - h(r),\]
and $g(0) = g(2) = 0$, $\displaystyle{\lim_{r \to \infty}}g(r) = \infty$, $h(0) = 0$, and if $r \in [0,2]$, then $g(r) \leq 0$, and $g(r) > 0$ in case of $r > 2$. Furthermore, $g(r)$ is monotone increasing on $r \in (\sqrt{2},\infty)$, and is monotone decreasing on $[0,\sqrt{2}]$. Furthermore, $h(r)$ is monotone (slowly) increasing on $r \in [0,\infty)$. By these observations on the drown graphs of such two functions $g(r)$  and $h(r)$, the graphs on $r > 0$ of $g(r)$ and $h(r)$ have only one common point $\delta > 2$, so that such only one point  $\delta > 2$ plays a role to satisfy the following property:  
%To do that we consider the function
%\begin{equation}\label{f-r}
 %f(r) = r^4 - 4r^2  -4 m^2 \log(1+r^{2\theta}), \quad  r \geq 0 .
%\end{equation}
%We notice that $f(0)=0$, 
%$$ \lim _{r \rightarrow \infty} (r^4 - 4r^2 -4 m^2 \log(1+r^{2\theta})) =+ \infty, $$
%for $ 0\leq \theta \leq 2 $,  and 
%$$ f'(r)= 4r^3 - 8r -4 m^2 \frac{2\theta r^{2\theta-1}}{ 1+r^{2\theta}}. $$ 
%\noindent
%t $f(r)$ has a global minimum in a unique  $r_0>0$  when  $ 0\leq \theta \leq 2 $. Hence,  one can conclude that there exists a unique $\delta=\delta(\theta) >0$ such that 
\begin{align}\label{delta}
 &f(r) = g(r)-h(r) = r^4 - 4r^2 -4 m^2 \log(1+r^{2\theta})) < 0, \quad \mbox{for} \;\; 0 < r <  \delta,\nonumber
\\
&f(r) = g(r)-h(r) = r^4 - 4r^2 -4m^2 \log(1+r^{2\theta})) \geq 0,  \quad \mbox{for} \;\;r \geq  \delta.
\end{align}
Note that thse arguments above hold good for each fixed $\theta \geq 0$, and the uniquely determined constant $\delta > 2$ depends on each fixed $\theta \geq 0$. Therefore, we conclude that the characteristics roots  are complex on the zone $\{|\xi| <\delta \}$ and are given by
$$ \lambda_{\pm}=a(\xi) \pm ib(\xi), \quad i=\sqrt{-1},$$
where
\begin{align} \label{a-b-low}
a(\xi) = -\frac{|\xi|^2}{2}, \quad b(\xi) = \frac{\sqrt{ 4 |\xi|^2 +4 m^2 \log(1+|\xi|^{2\theta}) -|\xi|^4} }{2},
\end{align}
and they are real on the region $\{|\xi| \geq\delta \} $, and can be defined by
$$\lambda_{\pm}=a(\xi) \pm d(\xi),$$
\begin{align} \label{a-d-low}
d(\xi) = \frac{\sqrt{\vert\xi\vert^{4} - 4 |\xi|^2 -4 m^2 \log(1+|\xi|^{2\theta})} }{2}.
\end{align}

%%%%%%%%%%%%%%%%
\subsection{Solution formula for low frequency and asymptotic profile}
At this point,  without loss of generality, we may assume that the initial data $u_0=0$. 

Then, the solution formula to the Cauchy problem \eqref{F-eq}--\eqref{F-ini} for low frequencies
\begin{align}
\hat{u}(t,\xi) = e^{-a(\xi) t} \frac{\sin (b(\xi) t) }{b(\xi)} \hat{u}_{1} + \frac{a(\xi)}{b(\xi)}e^{-a(\xi) t}\sin (b(\xi) t) \hat{u}_{0} + e^{-a(\xi) t} \cos (b(\xi) t) \hat{u}_{0}  \nonumber
\end{align}
becomes
%%%
\begin{align} \label{sol-formu}
\hat{u}(t,\xi) = e^{-a(\xi) t} \frac{\sin (b(\xi) t) }{b(\xi)}\hat{u}_{1}, \; t \geq 0,
\end{align}
where $a(\xi)$ and $b(\xi)$ are given  in \eqref{a-b-low} for $0\leq \theta$.

Our interest in this section is to find a candidate to be an asymptotic profile of the solution \eqref{sol-formu}. 

In that direction we consider the  elementary decomposition of the initial data $\hat{u}_1$ such that
%%%
$$\hat{u}_1(\xi)=A(\xi) - iB(\xi) +P_1,$$
where
\begin{equation}\label{A-B}
A(\xi)=\int_{{\bf R}^{n}}u_1(x)(\cos(\xi x) - 1)dx, \quad B(\xi)=\int_{{\bf R}^{n}}u_1(x)\sin(\xi x) dx, \quad P_1=\int_{{\bf R}^{n}}u_1(x)dx. 
\end{equation}
\noindent
Thus,
%%%
\begin{align*}
\hat{u}(t,\xi) = P_1e^{-a(\xi) t} \frac{\sin (b(\xi) t) }{b(\xi)} 
+\big(A(\xi) - iB(\xi)\big)e^{-a(\xi) t} \frac{\sin (b(\xi) t) }{b(\xi)} .
\end{align*}

Now, we apply the mean value theorem to write that
\begin{align*}
\sin (b(\xi) t)&=\sin\big(t\sqrt{|\xi|^2 +m^2\log(1+|\xi|^{2\theta} )}\big)\\
&+ t\big(b(\xi)-\sqrt{|\xi|^2+m^2 \log(1+|\xi|^{2\theta})}\big)\cos(\mu(\xi)t)
\end{align*}
with 
$\mu(\xi)=(1-\eta) b(\xi)+\eta\sqrt{|\xi|^2+m^2 \log(1+|\xi|^{2\theta})}$ 
for some $0<\eta<1$, $t>0$.
\vspace{0.3cm}
\noindent
The above expression for $\sin (b(\xi) t)$ implies that 
\begin{align*}\label{solu}
\hat{u}(t,\xi) &= P_1e^{-a(\xi) t} \frac{\sin\big(t\sqrt{|\xi|^2 +m^2\log(1+|\xi|^{2\theta} )} \big)}{b(\xi)} \nonumber\\
&+tP_1e^{-a(\xi)t}\frac{b(\xi)-\sqrt{|\xi|^2+m^2 \log(1+|\xi|^{2\theta})}}{b(\xi)}\cos(\mu(\xi)t)\\
&+\big(A(\xi) - iB(\xi)\big)e^{-a(\xi) t} \frac{\sin (b(\xi) t) }{b(\xi)}. 
\end{align*}

Now, in order to obtain a suitable expression of an asymptotic profile we observe that
\begin{equation*}\label{ere}
\frac{1}{b(\xi)} = \frac{1}{\sqrt{|\xi|^2+m^2 \log(1+|\xi|^{2\theta})}} +R(|\xi|),
\end{equation*}
where
\begin{equation}\label{ere}
 R(|\xi|) =\frac{|\xi|^4}{{2}b(\xi) \big(|\xi|^2+m^2 \log(1+|\xi|^{2\theta})\big) \left ( 2+ \sqrt{4- \frac{|\xi|^4}{|\xi|^2+m^2 \log(1+
 |\xi|^{2\theta})}} \right )}
\end{equation}
for $0<|\xi|<\delta$.

Then one can rewrite $\hat{u}(t,\xi)$ as 
%%%
%%%
\begin{align}\label{solu-for}
\hat{u}(t,\xi) = \varphi(t,\xi) +F_1(t,\xi) +F_2(t,\xi) +F_3(t,\xi), 
\end{align}
where 
%%%
\begin{align}\label{Func-Fi}
F_1(t,\xi)=& P_1e^{-a(\xi) t} R(\vert\xi\vert)\sin\big(t\sqrt{|\xi|^2 +m^2\log(1+|\xi|^{2\theta} )} \big), \nonumber
\\
F_2(t,\xi)=&tP_1e^{-a(\xi)t}\frac{b(\xi)-\sqrt{|\xi|^2+m^2 \log(1+|\xi|^{2\theta})}}{b(\xi)}\cos(\mu(\xi)t), \\
%%%
%%
F_3(t,\xi)=&\big(A(\xi) - iB(\xi)\big)e^{-a(\xi) t} \frac{\sin (b(\xi) t) }{b(\xi)},
\nonumber 
\end{align}
and the function $\varphi(t,\xi)$ is the asymptotic profile, which can be given by
%%%
\begin{equation}\label{Perfil}
\varphi(t,\xi)=P_1e^{-a(\xi) t} \frac{\sin\big(t\sqrt{|\xi|^2 +m^2\log(1+|\xi|^{2\theta} )} \big)}{\sqrt{|\xi|^2+m^2 \log(1+|\xi|^{2\theta})}}. 
\end{equation}

Next we need to estimate the $L^2$-norms of the remainder terms $F_i(t,\xi), i=1,2,3,$ on the zone of low frequency $\{|\xi|\leq\delta\}$.

Before beginning to estimate them, in order to avoid a singularity of $|\xi|^4/b(\xi)=|\xi|^4/b(|\xi|)$ in  $|\xi|=\delta$ we need to separate the zone of low frequency into two parts.
To do that we note the limit 
$$ \lim _{r \rightarrow +0} \frac{r^{\theta}
}{b(r)}= \lim _{r \rightarrow +0} \frac{r^{\theta}}{ \sqrt{4r^2+ 4 m^2 \log(1+r^{2\theta}) - r^4}}  < +\infty. $$
Thus one can take a small number $\delta_0, \;0<\delta_0<\delta$ such that 
\begin{equation}\label{delta-zero}
0\leq \frac{r^{\theta}}{b(r)}=\frac{r^{\theta}}{ \sqrt{4r^2+ 4 m^2 \log(1+r^{2\theta}) - r^4}}  \leq C 
\end{equation}
for all $|\xi|\leq \delta_0$, $0\leq \theta$, where $C$ is a positive constant depending only on $\theta$ and $\delta_0$ because $b(r)>0$ on the interval $(0,\delta_0]$. Therefore we need to consider seeral estimates on the low frequency zone $\{|\xi|\leq \delta_0\}$, and on the intermediate frequency zone, that is, the 
middle zone $\{\delta_0 \leq |\xi|\leq \delta\}$.

Now, we first consider $F_{1}$-estimates in the low frequency zone $|\xi| \leq \delta_0$.

\begin{align}\label{est-F1}
\int_{|\xi|\leq\delta_0} |F_1(t,\xi)|^2d\xi &\leq \int_{|\xi|\leq\delta_0} |P_1|^2e^{-2a(\xi) t} |R(\vert\xi\vert)|^2d\xi \nonumber\\
&=\omega_n P_1^2\int_0^{\delta_0}r^{n-1}e^{-tr^2}|R(r)|^2dr.
\end{align}

Here, we need to get an  estimate to $|R(r)|$. According to its definition in \eqref{ere} and the estimate \eqref{delta-zero} it follows that for $|\xi|\leq \delta_0$,
\begin{align}\label{R-estim}
|R(r)|&= \dfrac{r^4}{2\big|b(r)\big| \big(r^2+m^2 \log(1+r^{2\theta})\big) \left ( 2+ \sqrt{4- \dfrac{r^4}{r^2+m^2 \log(1+
 r^{2\theta})}} \right )} \nonumber\\
 &\leq 
 C\dfrac{r^{4-\theta}}{  \big(r^2+m^2 \log(1+r^{2\theta})\big) \left ( 2+ \sqrt{4- \dfrac{r^4}{r^2+m^2 \log(1+
 r^{2\theta})}} \right )} \nonumber\\
 &\leq 
 C\dfrac{r^{4-\theta}}{ 2 \big(r^2+m^2 \log(1+r^{2\theta})\big)  },
%%% 
\end{align}
because 
$ \sqrt{4- \dfrac{r^4}{r^2+m^2 \log(1+
 r^{2\theta})}} >0$ for $|\xi| \leq \delta_0$, where $C > 0$ is a generous constant (if necessarily).
 
 Combining the estimates \eqref{R-estim} and \eqref{est-F1} we arrive at
%%%
\begin{align}\label{est-F1c}
\int_{|\xi|\leq\delta_0} |F_1(t,\xi)|^2d\xi 
&=C^2\omega_n P_1^2\int_0^{\delta_0}r^{n-1}e^{-tr^2} \dfrac{r^{8-2\theta}}{ 4 \big(r^2+m^2 \log(1+r^{2\theta})\big)^2 }dr\nonumber\\
&\leq C^2 \omega_n P_1^2\int_0^{\delta_0}e^{-tr^2} \dfrac{r^{n+7-2\theta}
}{ 4 m^4 \log^2(1+r^{2\theta}) }dr\nonumber\\
&\leq C_{n,\delta_0,\theta}\frac{P_1^2}{m^{2}}\int_0^{\delta_0}e^{-tr^2} r^{n+7-6\theta}
dr\nonumber
\\
&\leq C_{n,\delta_0,\theta}\frac{P_1^2}{m^{2}} t^{-\frac{n+8-6\theta}{2}}, \quad t>0, 
\end{align} 
due to Lemma \ref{lem-t}, and the fact that the function $\dfrac{r^{4\theta}}{ \log^2(1+r^{2\theta}) }$ is bounded in the interval $(0,\delta_0]$ because $\displaystyle{\lim_{r \rightarrow +0}}\dfrac{r^{2\theta}}{ \log(1+r^{2\theta}) }$ is finite. The constant $C_{n,\delta_0,\theta}$ depends only on $n$ , $\delta_0$ and $\theta$.
\\
 
A next step is to estimate $F_2(t,\xi)$.
%%%5
\begin{align}\label{est-F2}
\int_{|\xi|\leq\delta_0} |F_2(t,\xi)|^2d\xi &= \int_{|\xi|<\delta_0} t^2P_1^2e^{-2a(\xi)t}\frac{|b(\xi)-\sqrt{|\xi|^2+m^2 \log(1+|\xi|^{2\theta})}|^2}{|b(\xi)|^2}\cos^2(\mu(\xi)t)d\xi \nonumber\\
&\leq t^2P_1^2\int_{|\xi|<\delta_0} e^{-t|\xi|^2}\frac{|b(\xi)-\sqrt{|\xi|^2+m^2 \log(1+|\xi|^{2\theta})}|^2}{|b(\xi)|^2}d\xi \nonumber\\
&\leq  \omega_n t^2P_1^2\int_0^{\delta_0} r^{n-1}e^{-tr^2}\frac{\big|b(r)-\sqrt{r^2+m^2 \log(1+r^{2\theta})}\big|^2}{|b(\xi)|^2}dr.
\end{align}
\noindent
To estimate the last integral above we note that
$$|b(\xi)-\sqrt{r^2+m^2 \log(1+r^{2\theta})}|=\frac{r^4}{4b(\xi)+4\sqrt{r^2+m^2 \log(1+r^{2\theta})}} \leq \frac{r^4}{\sqrt{r^2+m^2 \log(1+r^{2\theta})}}
$$
for $0<r\leq \delta_0$ because $b(\xi)$ is positive on this interval. Then,  as in \eqref{est-F1c} based on the estimate \eqref{delta-zero} we may conclude that
$$\frac{|b(r)-\sqrt{r^2+m^2 \log(1+r^{2\theta})}|}{|b(r)|} \leq 
\frac{Cr^{4-\theta}}{\sqrt{r^2+m^2 \log(1+r^{2\theta})}} \leq \frac{C}{m}\frac{r^{4-\theta}}{\sqrt{\log(1+r^{2\theta})}} \leq \frac{C}{m}r^{4-2\theta}
$$ 
for $0< |\xi| \leq \delta_0$, with $C>0$ depending only on $\delta_0$ and $\theta$. Combining the above estimate with \eqref{est-F2} one can obtain
\begin{align}\label{est-F2b}
\int_{|\xi|\leq\delta_0} |F_2(t,\xi)|^2d\xi 
&\leq \frac{C^2}{m^{2}}\omega_n t^2P_1^2\int_0^{\delta_0} r^{n-1}e^{-tr^2}r^{8-4\theta}dr\nonumber\\
&\leq  C_{n,\delta_0,\theta}\; \omega_n m^{-2}t^2P_1^2  t^{-\frac{n+8-4\theta }{2}}\nonumber\\
&\leq  C_{n,\delta_0,\theta}\; \omega_n m^{-2}P_1^2 \; t^{-\frac{n+4-4\theta }{2}}.
\end{align}

\begin{rem}
{\rm Due to the estimate \eqref{est-F2b}, to get a suitable decay rate for the function $F_2(t,\xi)$ it suffices to constrain $\theta$, at least, to the interval $0 < \theta \leq 1$ (see \cite{DI} for the case $\theta = 0$).}
\end{rem}

A step now is to estimate the function $F_3(t,\xi)$ on the low frequency zone $|\xi|\leq \delta_0$. Applying Lemma \ref{lema2.6} with $\kappa = \theta, \; 0\leq \theta \leq 1 $, we have 
\begin{align}\label{est-F3-A1}
\int_{|\xi|\leq\delta _0 } |F_3(t,\xi)|^2d\xi &= \int_{|\xi|<\delta_0}\Big|\Big(A(\xi) - iB(\xi)\Big)e^{-a(\xi) t} \frac{\sin (b(\xi) t) }{b(\xi)}\Big|^2 d\xi \nonumber\\
& \leq \int_{|\xi|\leq\delta_0 }\big|A(\xi) - iB(\xi)\big|^2 e^{-t|\xi|^2} \frac{1}{b(\xi)^2} d\xi  \nonumber \\
& \leq \int_{|\xi|\leq\delta}\big(|A(\xi)| + |B(\xi)|\big)^2 e^{-t|\xi|^2} \frac{1}{b(\xi)^2} 
d\xi  \\
& \leq (K+M)^2 \| u_1 \|_{1,\theta}^2 \int _{|\xi| \leq \delta _0} e^{-t|\xi|^2} \frac{|\xi|^{2\theta}}{b(\xi)^2} d\xi \nonumber\\
&=(K+M)^2 \omega _n \| u_1 \|_{1,\theta}^2 \int _{0 }^{\delta_0} e^{-t r^2} \frac{4 r^{2\theta}}{4 r^2 + 4m^2 \log (1+r^{2\theta})-r^4} r^{ n-1} dr .\nonumber
\end{align}
\noindent
Since $ 0 < \theta \leq 1 $, we have 
\begin{align*}
\lim_{r \rightarrow 0 } \frac{4r^{2\theta} }{4 r^2 + 4m^2 \log (1+r^{2\theta})-r^4} = \lim_{r \rightarrow 0 }\frac{4r^{2\theta} }{ r^{2\theta}   \left ( 4 r^{2-2\theta} + 4m^2 \frac{\log (1+r^{2\theta})}{r^{2\theta}}  -r^{4-2\theta} \right )} =  \frac{1}{m^2} . 
\end{align*}
\noindent
Then, due to $ \delta_0  < \delta  $, there exists a constant $ C>0  $ such that 
\begin{align}\label{Est-F3-A2}
 \Big | \frac{4r^{2\theta} }{4 r^2 + 4m^2 \log (1+r^{2\theta})-r^4} \Big | \leq \frac{C}{m^{2}}, \quad 0<r \leq \delta  .
\end{align}
Using estimate \eqref{Est-F3-A2} in  \eqref{est-F3-A1}, we obtain with the help of Lemma \ref{lem-t} the following estimate
\begin{align}\label{est-F3-A3}
\int_{|\xi|\leq\delta _0 } |F_3(t,\xi)|^2d\xi &\leq C(K+M)^2m^{-2} \omega _n \| u_1 \|_{1,\theta}^2 \int _{0 }^{1} e^{-t r^2}  r^{ n-1} dr  \nonumber \\
&\leq C_n m^{-2} \| u_1 \|_{1,\theta}^2\; t ^{-\frac{n}{2}}
\end{align}
with $C_n>0$ a constant depending on the  space dimension $n$. %The last inequality in  estimate \eqref{est-F3-A3} can be obtained using Lemma \ref{general-p} because $r^2 \approx\log{1+r^2}$.

\section{Estimate for middle frequency $  \delta _0 \leq |\xi| < \delta$ }
 In this zone, the solution is given by  
\begin{align*}
\hat{u}(t,\xi) = e^{-a(\xi) t} \frac{\sin (b(\xi) t) }{b(\xi)}\hat{u}_{1}, \; t \geq 0. 
\end{align*}

We remember that $b(\xi)=b(|\xi|)=0$ for $|\xi|=\delta$. But the singularity of $\hat{u}(t,\xi)$ in $|\xi|=\delta$ is removable. In fact, we know that $ |\sin a| \leq a $ for all $a \geq 0 $, thus one can get  
\begin{align}\label{mf-10}
\int _{\delta _0 \leq |\xi|< \delta } |\hat{u}(t,\xi)|^2 d\xi &=  \int _{\delta _0 \leq |\xi|< \delta } e^{-2a(\xi) t} \frac{\sin^2 (b(\xi) t) }{b(\xi)^2} |\hat{u}_{1}(\xi)|^2 d\xi \nonumber\\
&\leq t^2 \int _{\delta _0 \leq |\xi|\leq \delta } e^{-2 |\xi|^2 t} |\hat{u}_{1}(\xi)|^2 d\xi\nonumber\\
&\leq t^2 e^{-2\delta _0^2 t } \int _{\delta _0 \leq |\xi| \leq  \delta } |\hat{u}_{1}(\xi)|^2 d\xi \nonumber\\
&\leq t^2 e^{-2\delta _0^2\; t }  \|u_1\|^2, \quad t>0,
%%%%. 
\end{align} 
that is, the solution $\hat{u}(t,\xi)$  decays exponentially on the middle zone. This estimate holds for $ 0 \leq \theta \leq 2$.

%%%%%%XXXXXXXXXXXXXXXXXXXXXXXXXXXXXXXXXXXXX%%%%
\subsection{Estimate on high frequency $|\xi| \geq \delta$}
 In this zone the characteristics roots are real and given for $0\leq \theta$ by
 $$ \lambda _{\pm} = \frac{-|\xi|^2 \pm \sqrt{|\xi|^4 - 4 |\xi|^2 -4 m^2 \log(1+|\xi|^{2\theta}) }}{2},$$ 
 and the solution formula  in the high frequency region is
\begin{align}
\hat{u}(t,\xi) = e^{-a(\xi) t} \frac{\sinh(b(\xi) t) }{b(\xi)} \hat{u}_{1}, 
\end{align}
where 
\begin{align}
a(\xi) = \frac{|\xi|^2}{2}, \quad b(\xi) = \frac{\sqrt{ |\xi|^4 -4 |\xi|^2 - 4m^2 \log(1+|\xi|^{2\theta})  }}{2} .
\end{align}
\noindent
We recall that $\delta>0$ is such that $f(\delta)=0$  ($\delta >2$ is the unique positive zero of $f(r)$) and $f(r)$  is a increasing positive function  for $r>\delta$ defined    by
$$f(r)=r^4 - 4 r^2 -4 m^2 \log(1+r^{2\theta}).
$$
\noindent
Then there exists $\delta_1>\delta$ such that 
\begin{align} \label{Exi-delta-1}
0<\sqrt{1 - \frac{4}{\delta_1^2} -4 \frac{m^2}{\delta_1^4} \log(1+\delta_1^{2\theta})
 } <1, 
 \end{align}
 and, in particular, 
 \begin{align} \label{deltas}
\sqrt{1 - \frac{4}{\delta_1^2} -4 \frac{m^2}{\delta_1^4} \log(1+\delta_1^{2\theta})
 } \geq  
\sqrt{1 - \frac{4}{r^2} -4 \frac{m^2}{r^4} \log(1+r^{2\theta})
 }\geq 0
 \end{align}
 for $\delta \leq r\leq \delta_1.$
 
 Now we are able to estimate the $L^2$-norm of $\hat{u}$ on the region of high frequency.  We do that into two steps. The first step is to estimate on the middle frequency zone $\{\delta \leq \vert\xi\vert \leq \delta_{1}\}$:
 \begin{align*}
 \int_{\delta\leq |\xi\vert \leq \delta_1}|\hat{u}(t,\xi)|^2d\xi &\leq  \int_{\delta\leq |\xi| \leq \delta_1} |\hat{u}_{1}(t,\xi)|^2 e^{-2ta(\xi)}\frac{\sinh^2(b(\xi)t) }{b^2(\xi)}d\xi\nonumber\\
 &\leq  \|u_1\|_1^2\int_{\delta\leq |\xi| \leq \delta_1}  e^{-t|\xi|^2}\frac{\sinh^2(b(\xi)t) }{b^2(\xi)}d\xi\nonumber\\
  &\leq   \|u_1\|_1^2\int_{\delta\leq |\xi| \leq \delta_1} e^{-t|\xi|^2} t^2 e^{2b(\xi)t}d\xi\nonumber\\
  &\leq   \|u_1\|_1^2 \;t^2 \int_{\delta}^{\delta_1} \omega_n r^{n-1}e^{-tr^2}e^{2b(r)t}dr\nonumber\\
   &\leq \omega_n  \|u_1\|_1^2 \;t^2 \int_{\delta}^{\delta_1} r^{n-1}e^{-tr^2}e^{\sqrt{r^4 - 4 r^2 -4 m^2 \log(1+r^{2\theta})}t}dr\nonumber\\
    &\leq \omega_n \delta_1^{n-1}\|u_1\|_1^2 \;t^2 \int_{\delta}^{\delta_1} e^{-tr^2\big(1-{\sqrt{1 - \frac{4}{r^2} -4\frac{ m^2}{r^4} \log(1+r^{2\theta})}\big)}}dr,\nonumber\\
 \end{align*}
which holds for any space dimensions $n\geq 1$ and $t>0$. In the above estimates we have just used the elementary inequality Lemma \ref{lem-sen-hiper}.%$$\frac{\sinh(x)}{x}\leq e^x,  \mbox{for all} \; x>0.$$

By using \eqref{deltas} in  the last estimate above,  one can get
\begin{align}\label{sub-hi}
 \int_{\delta\leq |\xi \leq \delta_1}|\hat{u}(t,\xi)|^2d\xi  
  &\leq \omega_n \delta_1^{n-1}\|u_1\|_1^2 \;t^2 \int_{\delta}^{\delta_1} e^{-tr^2\big(1-{\sqrt{1 - \frac{4}{\delta_1^2} -4\frac{ m^2}{\delta_1^4} \log(1+\delta_{1}^{2\theta})}\big)}}dr\nonumber\\
  %%%
  &\leq \omega_n \delta_1^{n-1}\|u_1\|_1^2 \; t^2(\delta_1 - \delta) e^{-t\delta^2\big(1-{\sqrt{1 - \frac{4}{\delta_1^2} -4\frac{ m^2}{\delta_1^4} \log(1+\delta_{1}^{2\theta})}\big)}}\nonumber\\
  &\leq C_n\|u_1\|_1^2\; t^2 e^{-\alpha t}, \quad t>0,
 \end{align}
 %%%%%
where $\alpha=\delta^2\big(1-{\sqrt{1 - \frac{4}{\delta_1^2} -4\frac{ m^2}{\delta_1^4} \log(1+\delta_1^{2\theta})}}\big) > 0$ (due to \eqref{Exi-delta-1}),  and $C_n$ a positive constant depending only on $n$ and $\theta$ because $\delta$ depends on $\theta$.

Now, to complete the estimate on the zone of high frequency we proceed with the second step, that is, we estimate $\hat{u}(t, \cdot) $ on $|\xi|\geq \delta_1$. Then, because of  Lemma \ref{lem-sen-hiper} one can obtain
%%%%%
\begin{align}\label{hyper-hi}
 \int_{ |\xi \geq \delta_1}&|\hat{u}(t,\xi)|^2d\xi \leq  \int_{ |\xi| \geq \delta_1} |\hat{u}_{1}(t,\xi)|^2 e^{-2ta(\xi)}\frac{\sinh^2(b(\xi)t) }{b^2(\xi)}d\xi\nonumber\\
 &\leq  \|u_1\|_1^2\int_{ |\xi| \geq \delta_1}  e^{-2t|\xi|^2} \; t^2 e^{2b(\xi)t }d\xi\nonumber\\
 &\leq \omega_n \|u_1\|_1^2 \; t^2 \int_{ \delta_1}^{\infty} r^{n-1}e^{-tr^2}e^{t\sqrt{r^4 - 4 r^2 -4 m^2 \log(1+r^{2\theta})}}dr\nonumber\\
 &\leq \omega_n \|u_1\|_1^2 \; t^2 \int_{ \delta_1}^{\infty} r^{n-1}e^{-tr^2\big(1-\sqrt{1 - \frac{4} {r^2}-\frac{4 m^2}{r^4} \log(1+r^{2\theta})}\big)}dr\nonumber\\
  &\leq \omega_n \|u_1\|_1^2 \; t^2 \int_{ \delta_1}^{\infty} r^{n-1}e^{-\frac{tr^2}{2}\big(1-\sqrt{1 - \frac{4} {r^2}-\frac{4 m^2}{r^4} \log(1+r^{2\theta})}\big)}e^{-\frac{tr^2}{2}\big(1-\sqrt{1 - \frac{4} {r^2}-\frac{4 m^2}{r^4} \log(1+r^{2\theta})}\big)}dr
  \nonumber\\
  &\leq \omega_n \|u_1\|_1^2 \; t^2 \int_{ \delta_1}^{\infty} r^{n-1}e^{-\frac{r^2}{2}\big(1-\sqrt{1 - \frac{4} {r^2}-\frac{4 m^2}{r^4} \log(1+r^{2\theta})}\big)}e^{-\frac{{\bf t}\delta_1^2}{2}\big(1-\sqrt{1 - \frac{4} {\delta_{1}^2}-\frac{4 m^2}{\delta_{1}^4} \log(1+\delta_{1}^{2\theta})}\big)}dr \nonumber\\
 %%%%%
   &\leq \omega_n \|u_1\|_1^2 \;  t^2\Big(\int_{ \delta_1}^{\infty} r^{n-1}e^{-\frac{r^2}{2}\big(1-\sqrt{1 - \frac{4} {^2}-\frac{4 m^2}{r^4} \log(1+r^{2\theta})}}\big)dr\Big) e^{-\frac{{\bf t}\delta_1^2}{2}\big(1-\sqrt{1 - \frac{4} {\delta_1^2}-\frac{4 m^2}{\delta_1^4} \log(1+\delta_1^{2\theta})}\big)}\nonumber\\
   %%%%%
 &\leq \omega_n C_{n,\theta} \|u_1\|_1^2 \; t^2 e^{-\beta t},
\end{align}
which holds for $t \geq 1$. The constant 
 $\beta=\frac{\delta_1^2}{2}\big(1-{\sqrt{1 - \frac{4}{\delta_1^2} -4\frac{ m^2}{\delta_1^4} \log(1+\delta_1^{2\theta})}}\big)$ is a positive constant  (due to \eqref{Exi-delta-1}) depending only on $n$ and $\theta$ because $\delta_1$ depends on $\theta$.
% The constant $C_{n,\theta}>0$ is the value of the  last convergent integral from $\delta_1$ to $\infty$. 
%\vspace{0.6cm}
 
The estimates \eqref{sub-hi} and \eqref{hyper-hi}  prove the exponential decay of 
 the $L^2$-norm of $\hat{u}$ on the zone of high frequency.\\
 
%As consequences of a series of estimates \eqref{est-F1c}, \eqref{est-F2b}, \eqref{est-F3-A3}, \eqref{mf-10}, \eqref{sub-hi} and \eqref{hyper-hi} one can get the desired results stated in Theorem 1.1.

 %Finally we need to get an optimal decay estimate to the asymptotic profile $\varphi(t,\xi)$. In this direction we have
 %%
% \begin{pro} \label{prop-conclusion}
% There exists positive constant $C$  such that
 %\begin{align*}
 %% 
% \int |\hat{u}(t,\xi)- \varphi(t,\xi)|^2....d\xi \quad t>0.
 %%
% \end{align*}
 %\end{pro}
\section{Estimates to the asymptotic profile $\varphi(t,\xi)$   }

In this section one shall study the $L^{2}$-estimates of the profile $\varphi(t,x)$ itself in order to get sharp $L^{2}$-estimates of the solution itself. This is important to observe a singularity included in the solution $u(t,x)$. Such a singularity reflects on a ``growth rate'' as $t \to \infty$ of the quantity $\Vert u(t,\cdot)\Vert$ as stated in Theorem 1.2.  

\subsection{Decay estimates}

In this subsection we present several estimates on the asymptotic profile of the solution to problem \eqref{equation}--\eqref{initial}.

 To begin with, for low frequency parameter $ |\xi| \leq 1 $ one has  
 \begin{align*}
 \int _{|\xi|\leq \delta_0} |\varphi(t,\xi)|^2 d\xi &= P_1^2 \int _{0}^{\delta_0}  e^{-r^2t} \frac{\sin^2 (t\sqrt{r^2 + m^2 \log(1+r^{2\theta})} )}{r^2 + m^2 \log (1+r^{2\theta})} r^{n-1} dr\\
&\leq   P_1^2 \int _{0}^{1}  e^{-r^2t}  \frac{r^{2\theta}}{r^2 + m^2 \log (1+r^{2\theta})} r^{n-1-2\theta} dr,  
 \end{align*}
where $\delta_{0} > 0$ is a number already defined in previous section 4 (see \eqref{delta-zero}). Since $ 0 < \theta \leq 1 $, we have
\begin{align*}
 \lim _{r\rightarrow +0 } \frac{r^{2\theta}}{r^2 + m^2 \log (1+r^{2\theta})} = \lim _{r\rightarrow +0 }\frac{r^{2\theta}}{r^{2\theta} \left ( r^{2-2\theta} + m^2 \frac{\log (1+r^{2\theta})}{r^{2\theta} } \right ) } = \frac{1}{m^2} 
\end{align*}
so that there is a constant $C > 0$ such that 
$$ \frac{r^{2\theta}}{r^2 + m^2 \log (1+r^{2\theta})} \leq \frac{C}{m^{2}} , \quad 0 < r \leq 1.$$
Thus, we conclude that 
\begin{align*}
 \int _{|\xi|\leq \delta_0} |\varphi(t,\xi)|^2 d\xi &\leq \frac{C}{m^{2}}P_1^2 \int _{0}^{1}  e^{-r^2t} r^{n-1-2\theta} dr.
\end{align*}
\noindent
Applying Lemma \ref{lem-t} the next statement holds. 
%%5
\begin{lem}\label{lema-perfil} 
Let $0 < \theta \leq 1$  such that $n-2\theta>0$. Then,  it holds that
\begin{align}\label{perfil-low}
 \int _{|\xi|\leq \delta_0} |\varphi(t,\xi)|^2 d\xi 
\leq C_{n,\theta}\frac{P_1^2}{m^{2}} t^{-  \frac{n-2\theta}{2}}, \quad t\gg 1.
\end{align}
\end{lem}

\begin{rem}
{\rm The decay estimate \eqref{perfil-low} holds in the following cases:
%%%
\begin{itemize}
\item[{\rm (i)}] $n=1 $ with $ 0 < \theta < \frac{1}{2} $;  
\item[{\rm (ii)}] $n=2 $ with $ 0 < \theta < 1  $; 
\item[{\rm (iii)}] $ n \geq 3  $ and $0 < \theta \leq 1 $.
\end{itemize}
Note that the cases $n=1 $ with $\theta \in [\frac{1}{2},1]$ and and $ n=2$ with $\theta = 1$ are not considered in the classification above. These missing cases are studied  in the next Subsection 6.2.}
\end{rem}
  
  The effect on the decay estimates of the function $\varphi(t,\xi)$ in the zone of high frequency is exponentially small according to the next lemma.
  \begin{lem}\label{Lem-perfil-high}
  There exists a positive constant  $C_{\delta_0, n, \theta}$,   $0 \leq \theta \leq 1$, such that 
  \begin{align} \label{perfil-high}
  \int_{|\xi|\geq \delta_0}|\varphi(t,\xi)|^2 d\xi \leq C_{\delta_0, n,\theta} P_1^2\; e^{-\frac{\delta_0^2}{2}\;t} , \quad t \geq 1.
  \end{align}
  \end{lem}
{\it{Proof.}}  %Let  $C_{\delta_0}= \delta_0^2 +m^2 \log(1+\delta_0^{2\theta})$. Then
  \begin{align*}
   \int_{|\xi|\geq \delta_0}|\varphi(t,\xi)|^2 d\xi & 
   =\omega_nP^2_1 \int_{\delta_0}^{\infty} e^{-r^2t}r^{n-1}\frac{\sin^2\big({t\sqrt{r^2 + m^2 \log (1+r^{2\theta})}}\big)}{r^2 + m^2 \log (1+r^{2\theta})} dr\\
%%%   
   &\leq \omega_nP^2_1 \int_{\delta_0}^{\infty} e^{-r^2t}r^{n-1}\frac{1}{r^2 + m^2 \log (1+r^{2\theta})} dr\\
   &\leq \omega_nP^2_1 \int_{\delta_0}^{\infty} e^{-r^2t}r^{n-1}\frac{1}{\delta_0^{2}} dr\\
   &\leq \omega_n \frac{P^2_1}{\delta_{0}^{2}}\; e^{-\frac{\delta_0^2}{2}t} \int_{\delta_0}^{\infty} e^{-\frac{r^2}{2}}r^{n-1}dr, \quad t\geq 1.
  \end{align*}
  
This estimate proves the lemma because the  last integral is finite.
\hfill$\Box$
  %%%%
  \begin{rem}\label{Rem-perfil-above}
 {\rm The Lemmas \ref{lema-perfil} and \ref{perfil-high} imply that
  \begin{align}\label{perfil-above}
 \int _{{\bf R}^n} |\varphi(t,\xi)|^2 d\xi 
\leq C_{n,\theta,\delta_0}  \frac{P_1^2}{m^{2}}t^{-  \frac{n-2\theta}{2}}, \quad t>0,
\end{align}
with $C_{n,\theta,\delta_0}$ a positive constant independent of $t$.}
  \end{rem}

 The estimate obtained in \eqref{perfil-above} may be optimal, but to check it, one needs more precise estimates from below. The lemma below is one of main difficulties to be proved in the paper. 
  %%
%The  estimate in \eqref{perfil-above}  is  optimal but  to prove  that is true,  it rest to be showed  it from below.  
  %%
\begin{lem}\label{Lem-perfil-below}
For $0<\theta\leq 1$ with $n-2\theta>0$, there exist  positive constants  $C_{n,m, \theta}$  and $t_0$ such that 
\begin{align}\label{Lem-perfil-low}
\int _{{\bf R}^n} |\varphi(t,\xi)|^2 d\xi 
\geq C_{n,m,\theta}  P_1^2  t^{-  \frac{n-2\theta}{{2}}}, \quad t \geq t_0.
\end{align}
\end{lem}
{\it{Proof.}}
In order to capture the precise estimate from below we need to estimate the $L^2$-norm of the profile  only on the low frequency zone due to that  on the high frequency the $L^2$-norm of the profile decays exponential, according to Lemma \eqref{Lem-perfil-high}.

We first consider the change of variable $ w = \sqrt{r^2 +m^2\log (1+r^{2\theta})} $. Then $$ dw = \frac{r + \theta m^2 \frac{r^{2\theta -1}}{1+r^{2\theta}}}{\sqrt{r^2 +m^2\log (1+r^{2\theta})}} dr,$$
and
\begin{align*}
\int_{|\xi| \leq 1} &|\varphi(t,\xi)|^2 d\xi =  \omega_n P_1^2 \int _{0}^{1}  e^{-r^2t} \frac{\sin^2 (t\sqrt{r^2 + m^2 \log(1+r^{2\theta})} )}{r^2 + m^2 \log (1+r^{2\theta})} r^{n-1}  dr \\
&= \omega_n P_1^2 \int _{0}^{1}   \frac{ e^{-r^2t} r^{n-1} \sin^2 (t\sqrt{r^2 + m^2 \log(1+r^{2\theta})} )}{\sqrt{r^2 + m^2 \log (1+r^{2\theta})} \Big ( r + \theta m^2 \frac{r^{2\theta -1}}{1+r^{2\theta}} \Big )}   \frac{r + \theta m^2 \frac{r^{2\theta -1}}{1+r^{2\theta}}}{\sqrt{r^2 +m^2\log (1+r^{2\theta})}} dr . 
\end{align*}  
For  $ 0 \leq r \leq 1 $ and $ 0 \leq \theta \leq 1 $, we have 
\begin{align}\label{A2_03maio}
\frac{m^2}{2} r^{2\theta} \leq w^2=  r^2 +m^2\log (1+r^{2\theta}) \leq (1+m^2) r^{2\theta}, 
\end{align}
%%%
that implies  the important equivalence  $ r \approx w^{\frac{1}{\theta}} $, for $ 0 \leq r \leq 1 $. We also notice that
$$ r + \theta m^2 \frac{r^{2\theta -1}}{1+r^{2\theta}} \approx r^{2\theta -1}, \quad  0 < r \leq 1, \quad 0 \leq \theta \leq 1 .$$
Thus, the above change of variable give us the estimate
\begin{align*}
\int_{|\xi| \leq 1} |\varphi(t,\xi)|^2 d\xi  &\approx \omega_{n}P_1^2 \int _{0}^{\eta } e^{- w^{\frac{2}{\theta}} t} \frac{\sin^2 (t w)}{w } w^{\frac{n-2\theta}{\theta}} dw \\
&=  \omega_{n}P_1^2 \int _{0}^{\eta } e^{- w^{\frac{2}{\theta}} t} \sin^2 (tw) w^{\frac{n-2\theta}{\theta}-1} dw \\
&=  \omega_{n}P_1^2 \int _{0}^{\eta } e^{- w^{\frac{2}{\theta}} t} \sin^2 \Big ( t^{1-\frac{\theta}{2}} t^{\frac{\theta}{2}}w \Big ) w^{\frac{n-3\theta}{\theta}} dw,
\end{align*}
where $ \eta = \sqrt{1+ m^2 \log 2} $. 

Now we set $ y = t^{\frac{\theta}{2}}w  $. Then for $t \gg 1$ %$ w = t^{-\frac{\theta}{2}} y $, and
\begin{align*}
\int_{|\xi| \leq 1} |\varphi(t,\xi)|^2 d\xi  &\approx \omega_{n}P_1^2 t^{-\frac{\theta}{2}} \int _{0}^{\eta } e^{- w^{\frac{2}{\theta}} t} \sin^2 \Big ( t^{1-\frac{\theta}{2}} t^{\frac{\theta}{2}}w \Big ) w^{\frac{n-3\theta}{\theta}} t^{\frac{\theta}{2}} dw\\ 
& = \omega_{n}P_1^2  t^{-\frac{n-3\theta }{2}} t^{-\frac{\theta}{2}}  \int _{0}^{\eta t^{\frac{\theta}{2}}} e^{-y^{\frac{2}{\theta}}} \sin^2 \Big ( t^{1-\frac{\theta}{2}}  y   \Big ) y ^{\frac{n-3\theta}{\theta}} dy \\
&\geq \omega_{n}P_1^2  t^{-\frac{n-2\theta }{2}}  \int _{0}^{1} e^{-y^{\frac{2}{\theta}}} \sin^2 \Big ( t^{1-\frac{\theta}{2}}  y   \Big ) y ^{\frac{n-3\theta}{\theta}} dy \\
&= \omega_{n}\frac{P_1^2}{2}  t^{-\frac{n-2\theta }{2}} \left ( \int _{0}^{1} e^{-y^{\frac{2}{\theta}}}  y ^{\frac{n-3\theta}{\theta}} dy  - \int _{0}^{1} e^{-y^{\frac{2}{\theta}}} y ^{\frac{n-3\theta}{\theta}} \cos \Big ( 2t^{1-\frac{\theta}{2}}  y   \Big )  dy \right )\\
&=: \omega_{n}\frac{P_1^2}{2}  t^{-\frac{n-2\theta }{2}} (A_{n} - F_{n}(t)),
\end{align*}
where 
\begin{align*}
&A_n:= \int _{0}^{1} e^{-y^{\frac{2}{\theta}}}  y ^{\frac{n-3\theta}{\theta}} dy, \\
&F_n(t) := \int _{0}^{1} e^{-y^{\frac{2}{\theta}}} y ^{\frac{n-3\theta}{\theta}} \cos \Big ( 2t^{1-\frac{\theta}{2}}  y   \Big )  dy.
\end{align*}
\noindent
For $ 0 <  \theta \leq 1  $ and $ n-2\theta > 0 $ the function  $  y \mapsto e^{-y^{\frac{2}{\theta}}} y ^{\frac{n-3\theta}{\theta}}$ is integrable on the interval  $[0,1]$. Thus, it follows from the Riemann-Lebesgue theorem that
$$ \lim _{t \rightarrow \infty } F_n(t) =0 , $$
due to $ t^{1-\frac{\theta}{2}} \rightarrow \infty $ as $ t \rightarrow \infty $. Therefore, there exists $ T>0 $ such that 
$$ -\frac{A_n}{2} \leq F_n(t) \leq \frac{A_n}{2}, \quad t \geq T,$$ 
so that  
$$ A_n -  F_n(t) \geq \frac{A_n}{2}  , \quad t \geq T. $$ 
Therefore
\begin{align*}
\int_{|\xi| \leq 1} |\varphi(t,\xi)|^2 d\xi &\geq \omega_{n}\frac{A_n P_1^2}{4}  t^{-\frac{n-2\theta }{2}}, \quad t \geq T,
\end{align*}
which implies the desired statement.  
\hfill $\Box $

%%%%%%%%%%%%%%%%%%%%%%%%%%%%%%%%%%%%%%%%%%%%%%%%%%%%%%%%%%%%%%%%%%%%%%%%%%%%%%%%%%%%%%%%%%%%%%%%%%%%%%%%%%%%%%%%%%%%%%%%

\subsection{Infinite time blow-up estimates}\label{blowup}
The omitted cases $n=1 $ with $ \frac{1}{2} \leq \theta \leq 1$ and $ n=2 $ with $\theta =1 $ are more delicate and we discuss each of them in this subsection through some lemmas. We prove that the asymptotic profile of the solution blows up in infinite time. This is a quite novel property which reflects a singularity hidden in the solution itself.

\subsubsection*{(I)\,Estimates for $n=1$ and $\theta \in (1/2,1]$.}

To obtain sharp estimates to this cases, we first prepare the following fact.
\begin{lem}\label{lema_erf} Let $c > 0$ and $\alpha > 0 $ be constants. Then it holds that
$$ \rm{erf}(c t^{-\alpha }) \sim t^{-\alpha }, \quad  t \geq c^{\frac{1}{\alpha}},$$
where $\rm{erf}(\cdot) $ is the \textit{Gauss error function} defined by
$$ {\rm erf} (x) = \frac{2}{\sqrt{\pi}} \int _{0}^{x} e^{-w^2 }dw.$$
\end{lem}
\textit{Proof.} Using the change of variable $ v= \frac{t^{\alpha}w}{c} $, we  have 
\begin{align}\label{A1_03maio}
\text{erf}(ct^{-\alpha }) &= \frac{2}{\sqrt{\pi}} \int _{0}^{ct^{-\alpha}} e^{-w^2 }dw = \frac{2c}{\sqrt{\pi} t^{\alpha }} \int_{0}^{1} e^{- \frac{c^2v^2}{t^{2\alpha}} }dv . 
\end{align}
Then, we can get the following upper  estimate
\begin{align*}
\text{erf}(ct^{-\alpha }) &\leq \frac{2c}{\sqrt{\pi} t^{\alpha }} \int_{0}^{1} dv = \frac{2c}{\sqrt{\pi}} t^{- \alpha }
\end{align*}
for $t>0$.
On the other hand, from \eqref{A1_03maio} and assuming that $ t \geq c^{\frac{1}{\alpha}} $, we have the following lower estimate
\begin{equation}\label{erf-10}
\text{erf}(ct^{-\alpha }) \geq \frac{2c}{\sqrt{\pi} t^{\alpha }} \int_{0}^{1} e^{- v^2 } dv = c\;{\rm erf}(1)t^{-\alpha }. 
\end{equation}
This completes the desired estimates.
\hfill $\Box$
\vspace{0.2cm}

In order to prove the following lemma we observe that for $n=1$ 
\begin{align}
\int _{|\xi|\leq 1} |\varphi(t,\xi)|^2 d\xi &= P_1^2 \int _{0}^{1}  e^{-r^2t} \frac{\sin^2 (t\sqrt{r^2 + m^2 \log(1+r^{2\theta})} )}{r^2 + m^2 \log (1+r^{2\theta})}  dr .
\end{align}
In the proof of the following lemma, we use the fact that 
\begin{align}\label{A2_03maio}
 \frac{m^2}{2} r^{2\theta} \leq  r^2 +m^2\log (1+r^{2\theta}) \leq (1+m^2) r^{2\theta} , 
\end{align}
due to 
\begin{align*}
&\frac{1}{2} r^{2\theta} \leq \log (1+r^{2\theta}) \leq r^{2\theta}, \\ 
&r^2 \leq r^{2\theta}
\end{align*}
for $ 0\leq r \leq 1  $ and $ 0\leq \theta \leq 1 $.

\begin{lem}\label{Lema6.3}
Let $ n=1$ and $ \frac{1}{2}< \theta\leq 1  $. There exist positive constants $C_1, C_2$ such that 
$$ \frac{C_1}{(1+m^{2})^{\frac{1}{2\theta}}} P_1^2 t^{\frac{2\theta -1 }{\theta}} \leq  \int _{|\xi|\leq 1} |\varphi(t,\xi)|^2 d\xi \leq \frac{C_{2}}{m^{2}}\frac{1}{2\theta -1 } t^{\frac{2\theta -1 }{\theta}} , \quad t\gg 1.$$
\end{lem}
\textit{Proof.}
First we get the upper estimate. 
\begin{align*}
\int _{|\xi|\leq 1} |\varphi(t,\xi)|^2 d\xi &= P_1^2 \int _{0}^{t^{-\frac{1}{\theta}}}  e^{-r^2t} \frac{\sin ^2 (t\sqrt{r^2 + m^2 \log(1+r^{2\theta})} )}{r^2 + m^2 \log (1+r^{2\theta})}  dr \nonumber \\
&+ P_1^2 \int _{t^{-\frac{1}{\theta}}}^{1}  e^{-r^2t} \frac{\sin ^2 (t\sqrt{r^2 + m^2 \log(1+r^{2\theta})} )}{r^2 + m^2 \log (1+r^{2\theta})}  dr \\
&=: P_1^2 \Big ( K_1(t) + K_2(t) \Big ).
\end{align*}

Set $L := \displaystyle{\sup_{\nu \ne 0}\vert\frac{\sin\nu}{\nu}\vert} < +\infty$. Then, one can estimate as follows: 
\begin{align}
K_1(t) &\leq L^{2}t^2 \int _{0}^{t^{-\frac{1}{\theta}}}  e^{-r^2t}   dr \leq L^2 t^2 \int _{0}^{t^{-\frac{1}{\theta}}}    dr \nonumber \\
&\leq L^{2}t^{2}t^{-\frac{1}{\theta}} \nonumber \\
&= L^{2}t^{\frac{2\theta -1}{\theta}}. \label{A3_03maio}
\end{align}

Due to $ |\sin a| \leq 1 $,  using \eqref{A2_03maio} and integration by parts, assuming that $ \theta \neq \frac{1}{2} $ one can get an upper bound estimate to $K_2(t)$  as follows: 
\begin{align}\label{A4_03maio}
K_2(t) &= \int _{t^{-\frac{1}{\theta}}}^{1}  e^{-r^2t} \frac{\sin ^2 (t\sqrt{r^2 + m^2 \log(1+r^{2\theta})} )}{r^2 + m^2 \log (1+r^{2\theta})}  dr  \nonumber \\
&\leq \int _{t^{-\frac{1}{\theta}}}^{1}  e^{-r^2t} \frac{1}{r^2 + m^2 \log (1+r^{2\theta})}  dr \nonumber  \\
&\leq \frac{2}{m^{2}}\int _{t^{-\frac{1}{\theta}}}^{1}\frac{1}{r^{2\theta}} dr \nonumber \\
&=\frac{2}{m^{2}}\left(\frac{1}{2\theta-1}t^{\frac{2\theta-1}{\theta}}-\frac{1}{2\theta-1}\right) \nonumber \\
&\leq \frac{2}{m^{2}}\frac{1}{2\theta-1}t^{\frac{2\theta-1}{\theta}}.
%&=\frac{2}{m^{2}}\int _{t^{-\frac{1}{\theta}}}^{1} (1+r^{2})^{-t}  \frac{1}{r^{2\theta}} dr \nonumber \\
%&= \frac{2}{m^{2}}\left(\frac{1}{1-2\theta }(1+r^{2})^{-t} r^{-2\theta +1} \Big |_{t^{-\frac{1}{\theta}}}^{1} - \frac{2t}{2\theta -1} \int _{t^{-\frac{1}{\theta}}}^{1} (1+r^2)^{-t-1}  r^{2-2\theta }  dr\right) \nonumber \\
%&\leq \frac{2}{m^{2}} \frac{1}{2\theta -1 }  \Big(1+ t^{-\frac{2}{\theta}} \Big )^{-t} t^{\frac{2\theta-1}{\theta }}. 
\end{align}
%\begin{align}
%K_2(t) &= \int _{t^{-\frac{1}{\theta}}}^{1}  e^{-r^2t} \frac{\sin ^2 (t\sqrt{r^2 + m^2 \log(1+r^{2\theta})} )}{r^2 + m^2 \log (1+r^{2\theta})}  dr  \nonumber \\
%&\leq \int _{t^{-\frac{1}{\theta}}}^{1}  e^{-r^2t} \frac{1}{r^2 + m^2 \log (1+r^{2\theta})}  dr \nonumber  \\
%&\leq \frac{2}{m^{2}}\int _{t^{-\frac{1}{\theta}}}^{1}   \frac{e^{-t\log (1+r^2)}}{r^{2\theta}} dr \nonumber \\
%&=\frac{2}{m^{2}}\int _{t^{-\frac{1}{\theta}}}^{1} (1+r^{2})^{-t}  \frac{1}{r^{2\theta}} dr \nonumber \\
%&= \frac{2}{m^{2}}\left(\frac{1}{1-2\theta }(1+r^{2})^{-t} r^{-2\theta +1} \Big |_{t^{-\frac{1}{\theta}}}^{1} - \frac{2t}{2\theta -1} \int _{t^{-\frac{1}{\theta}}}^{1} (1+r^2)^{-t-1}  r^{2-2\theta }  dr\right) \nonumber \\
%&\leq \frac{2}{m^{2}} \frac{1}{2\theta -1 }  \Big(1+ t^{-\frac{2}{\theta}} \Big )^{-t} t^{\frac{2\theta-1}{\theta }}. 
%\end{align}
From the estimates above one can observe that the Gauss kernel $e^{-tr^{2}}$ is no effective on the growth rates.

%For $ \theta > 1/2 $, it holds that $  \Big(1+ t^{-\frac{2}{\theta}} \Big )^{-t} \rightarrow 1 $  when $ t\rightarrow \infty  $. Then 
%\begin{align}
%K_2(t) \leq C \frac{1}{2\theta -1 } t^{\frac{2\theta -1 }{\theta }}, \quad t \gg 1 \label{A4_03maio}
%\end{align}
From \eqref{A3_03maio} and \eqref{A4_03maio}, we prove that if $\theta >\frac{1}{2} $ there exists a positive constant $C$ such that 
$$ \int _{|\xi|\leq 1} |\varphi(t,\xi)|^2 d\xi \leq P_{1}^{2}\frac{C}{m^{2}}\frac{1}{2\theta -1 } t^{\frac{2\theta -1 }{\theta }}, \quad t\gg 1 . $$ 
\vspace{0.2cm}

%For $n=1$ we have
%\begin{align}
%\int _{|\xi|\leq 1} |\varphi(t,\xi)|^2 d\xi &= P_1^2 \int _{0}^{1}  e^{-r^2t} \frac{\sin (t\sqrt{r^2 + %m^2 \log(1+r^{2\theta})} )}{r^2 + m^2 \log (1+r^{2\theta})}  dr .
%\end{align}

%For $0\leq r \leq 1 $ and $\theta \geq 1/2$ we may see that 
%$$ \frac{1}{2} r^{2\theta } \leq \log (1+r^{2\theta}) \leq  r^{2\theta} . $$
%Then  $ \frac{1}{2} \leq \theta \leq 1  $ and $0 \leq  r \leq 1$, we have 
%\begin{align}
%r^{2}+ m^2 \log (1+r^{2\theta}) \leq r^{2}+m^2 r^{2\theta} \leq (1+m^2) r^{2\theta} 
%\end{align}
%\begin{align}
%t \sqrt{r^{2}+ m^2 \log (1+r^{2\theta}) } \leq tr^{\theta}\sqrt{1+m^2}  . 
%\end{align}
In order to obtain the lower estimate we define
$$ r_0(t):= \frac{2^{\frac{1}{\theta}}}{(1+m^2)^\frac{1}{2\theta} } t^{-\frac{1}{\theta}}.$$
Then for  $ 0 \leq   r \leq  r_0(t) < 1$ we have
\begin{align}
 t \sqrt{r^{2}+ m^2 \log (1+r^{2\theta}) } \leq t\sqrt{1+m^2} r^{\theta} \leq 2, \quad t \gg 1,
\end{align}
due to \eqref{A2_03maio}. We observe that there exist a constant $C>0$ such that 
$ C \leq \displaystyle{\frac{\sin^{2}a}{a^2}} \leq L^{2}  $  for $0 < |a| \leq 2$. Thus for $ 0 \leq   r \leq  r_0(t)  $
\begin{align}
  C \leq \frac{\sin^2 (t \sqrt{r^{2}+ m^2 \log (1+r^{2\theta}) } ) }{t^2 (r^{2}+ m^2 \log (1+r^{2\theta}) ) }\leq L^{2}. 
\end{align}
Therefore, 
\begin{align}\label{a26junho1}
\int _{|\xi|\leq 1} |\varphi(t,\xi)|^2 d\xi &= P_1^2 \int _{0}^{1}  e^{-r^2t} \frac{\sin^{2}(t\sqrt{r^2 + m^2 \log(1+r^{2\theta})} )}{r^2 + m^2 \log (1+r^{2\theta})}  dr \nonumber \\
&\geq P_1^2 \int _{0}^{r_0(t) }  e^{-r^2t} \frac{\sin^{2}(t\sqrt{r^2 + m^2 \log(1+r^{2\theta})} )}{r^2 + m^2 \log (1+r^{2\theta})}  dr \nonumber \\ 
& \geq C t^2 P_1^2 \int _{0}^{r_0(t) }  e^{-r^2t} dr \nonumber \\
& = Ct^{\frac{3}{2}}  P_1^2 \int _{0}^{\sqrt{t} r_0(t) }  e^{-w^2} dw \nonumber \\
& = C P_1^2  t^{\frac{3}{2}}  \text{erf}\left ( \frac{2^{\frac{1}{\theta}}}{(1+m^2)^\frac{1}{2\theta} } t^{ \frac{1}{2} -\frac{1}{\theta}} \right ) \nonumber \\
%&\sim  P_1^2  t^{\frac{3}{2}} t^{ \frac{1}{2} -\frac{1}{\theta}}, \quad t\gg 1 \nonumber \\
&\geq  K\frac{P_1^2}{(1+m^{2})^{\frac{1}{2\theta}}} t^{\frac{2\theta -1}{\theta}}, \quad t\gg 1,
\end{align}
with some generous constant $K > 0$, where we just used the estimate \eqref{erf-10}. Note that \eqref{a26junho1} makes sense even if one chooses $m = 0$.
\hfill $\Box$

\subsubsection*{(II)\,Estimate for $n=1$ and $ \theta = 1/2 $.}
In this part we shall study the critical case of $\theta = 1/2$ just omitted in the statement of Lemma \ref{Lema6.3}. For this, we prepare the following lemma.  

\begin{lem}\label{LemmaA1} There exists a positive constant $ C $ such that
\begin{align}
\int _0^\infty \frac{e^{-w^2  } \sin^2(\sqrt{t} w)}{w} dw \geq C \log t, \quad t \gg 1. 
\end{align}
\end{lem}
\textit{Proof.}\,The whole idea comes from \cite{IO}. \\
Let $$ \nu _j = \left ( \frac{1}{4}+j \right )\frac{\pi}{\sqrt{t}}, \qquad \nu ' _j =\left ( \frac{3}{4}+j \right )\frac{\pi}{\sqrt{t}}, \quad j=0, 1, 2, 3, \cdots $$
and observe that 
$$ \sin ^2 (\sqrt{t} w) \geq  \frac{1}{2}, \forall w \in [\nu _j, \nu '_j].$$
Then 
\begin{align}\label{A5}
\int _0^\infty \frac{e^{-w^2} \sin^2(\sqrt{t} w)}{w} dw &\geq \sum _{j=0}^{\infty } \int _{\nu _j}^{\nu '_j} \frac{e^{-w^2} \sin^2(\sqrt{t}w)}{w} dw \nonumber \\
&\geq \frac{1}{2} \sum _{j=0}^{\infty } \int _{\nu _j}^{\nu '_j} \frac{e^{-w^2  } }{w} dw.
\end{align}
The function $ (0,\infty) \ni w \mapsto \displaystyle{\frac{e^{-w^2  } }{w}}$ is decreasing and the intervals $ [\nu _j, \nu' _j] $ and $ [ \nu '_j, \nu _{j+1} ] $ have the same length, then 
\begin{align}\label{A6}
 2 \sum _{j=0}^{\infty } \int _{\nu _j}^{\nu '_j} \frac{e^{-w^2  } }{w} dw \geq \sum _{j=0}^{\infty } \int _{\nu _j}^{\nu '_j} \frac{e^{-w^2  } }{w} dw + \sum _{j=0}^{\infty } \int _{\nu ' _j}^{\nu _{j+1}} \frac{e^{-w^2  } }{w} dw = \int _{\frac{\pi}{4 \sqrt{t}}}^{\infty} \frac{e^{-w^2  } }{w} dw.  
\end{align}
From \eqref{A5} and \eqref{A6}, for large $t > 0$ it follows that 
\begin{align}
\int _0^\infty \frac{e^{-w^2  } \sin^2(\sqrt{t}w)}{w} dw &\geq \frac{1}{4} \int_{\frac{\pi}{4\sqrt{t}} }^\infty \frac{e^{-w^2  } }{w} dw \nonumber \\
&\geq \frac{1}{4} \int _{\frac{\pi}{4\sqrt{t}} }^1 \frac{e^{-w^2  } }{w} dw \nonumber \\
&\geq \frac{1}{4}e^{-1}\int_{\frac{\pi}{4\sqrt{t}}}^{1}\frac{1}{w}dw \nonumber \\
&= \frac{1}{4}e^{-1}\left(\frac{1}{2}\log t - \log\frac{\pi}{4} \right).
%&= e^{-w^2} \log w \Big |_{\frac{\pi}{4\sqrt{t}} }^\infty + \frac{1}{2} \int _{\frac{\pi}{4\sqrt{t}} }^\infty w \log (w) e^{-w^2} dw \nonumber \\
%&= - e^{-\frac{\pi ^2}{16 t^2}} \log \Big (  \frac{\pi}{4\sqrt{t}}\Big ) + \frac{1}{2} \int _{\frac{\pi}{4\sqrt{t}} }^\infty w \log (w) e^{-w^2} dw .
\end{align}
%\begin{align}
%\int _0^\infty \frac{e^{-w^2  } \sin^2(\sqrt{t}w)}{w} dw &\geq \frac{1}{4} \int _{\frac{\pi}{4\sqrt{t}} }^\infty \frac{e^{-w^2  } }{w} dw \nonumber \\
%&\geq \frac{1}{4} \int _{\frac{\pi}{4\sqrt{t}} }^1 \frac{e^{-w^2  } }{w} dw \nonumber \\
%&= e^{-w^2} \log w \Big |_{\frac{\pi}{4\sqrt{t}} }^\infty + \frac{1}{2} \int _{\frac{\pi}{4\sqrt{t}} }^\infty w \log (w) e^{-w^2} dw \nonumber \\
%&= - e^{-\frac{\pi ^2}{16 t^2}} \log \Big (  \frac{\pi}{4\sqrt{t}}\Big ) + \frac{1}{2} \int _{\frac{\pi}{4\sqrt{t}} }^\infty w \log (w) e^{-w^2} dw .
%\end{align}
%But the integral  
%$$ \int _{0}^\infty w \log (w) e^{-w^2} dw $$
%is convergent. Then 
%\begin{align}
%\int _0^\infty \frac{e^{-w^2  } \sin^2(t w)}{w} dw &\geq - e^{-\frac{\pi ^2}{16 t}} \log \Big (  \frac{\pi}{4\sqrt{t}}\Big ) \\
%& \frac{1}{2} e^{-\frac{\pi ^2}{16 t}} \log (t) -  e^{-\frac{\pi ^2}{16 t}} \log \Big ( \frac{\pi}{4} \Big) .
%\end{align}
%Since 
%$$ \lim _{t \rightarrow \infty  } e^{-\frac{\pi ^2}{16 t}}  =1,$$
%there exists $T>0 $ such that 
%$$ \int _0^\infty \frac{e^{-w^2  } \sin^2(t w)}{w} dw \geq \frac{1}{4} \log (t) , \quad t \geq T .$$ 
This implies the desired estimate.
\hfill $\Box$

The next lemma state the $\log$-growth rate of the critical case for $\theta = 1/2$ and $n = 1$. 

\begin{lem}\label{Lema6.4}
Let $ n=1$ and $  \theta = \frac{1}{2}  $. There exist positive constants $C_1, C_2$ such that 
$$ \frac{C_1}{m^{2}+2}P_1^2 \log t  \leq  \int _{|\xi|\leq 1} |\varphi(t,\xi)|^2 d\xi \leq \frac{C_2}{m^{2}}P_1^2 \log t , \quad t\gg 1.$$
%where
%\[C^{*} := 2\left(\max_{0\leq r\leq1}\{\frac{m^{2}}{r+1}+2r\}\right)^{-1}.\]
\end{lem}
\textit{Proof.}
For $n=1$ and $ \theta = 1/2 $ and  $ t >1$, we have
\begin{align}
\int _{|\xi|\leq 1} |\varphi(t,\xi)|^2 d\xi &= P_1^2 \int _{0}^{1}  e^{-r^2t} \frac{\sin^2 (t\sqrt{r^2 + m^2 \log(1+r^{2\theta})} )}{r^2 + m^2 \log (1+r^{2\theta})}  dr \nonumber \\
&= P_1^2 \int _{0}^{\frac{1}{t^2}}  e^{-r^2t} \frac{\sin^2 (t\sqrt{r^2 + m^2 \log(1+r)} )}{r^2 + m^2 \log (1+r)}  dr\nonumber \\
&= P_1^2 \int _{0}^{\frac{1}{t^2}}  e^{-r^2t} \frac{\sin^2 (t\sqrt{r^2 + m^2 \log(1+r)} )}{r^2 + m^2 \log (1+r)}  dr\nonumber \\
 &+ P_1^2 \int _{\frac{1}{t^2}}^{1}  e^{-r^2t} \frac{\sin^2 (t\sqrt{r^2 + m^2 \log(1+r)} )}{r^2 + m^2 \log (1+r)}  dr \nonumber \\
&=: K_1(t) + K_2(t). \label{A1}
\end{align}
One first considers the upper bound estimates of \eqref{A1}.\\
Due to $ \sin^{2}a \leq a^{2}  $ for $ a \geq 0 $, one has
\begin{align}
K_1(t) = P_1^2 \int _{0}^{\frac{1}{t^2}}  e^{-r^2t} \frac{\sin^2 (t\sqrt{r^2 + m^2 \log(1+r)} )}{r^2 + m^2 \log (1+r)}  dr 
\leq  t^2 P_1^2 \int _{0}^{\frac{1}{t^2}}   dr  = P_1^2. \label{A2}
\end{align}
\noindent
In order to get an upper estimate to $K_2(t)$, we again remember that for $ 0\leq r \leq 1$,
\begin{align}\label{A4}
\frac{m^2}{2} r \leq r^2 + m^2 \log (1+r)  \leq (1+m^2) r,
\end{align}
as just observed in \eqref{A2_03maio} with $2\theta = 1$. Therefore
\begin{align}
K_2(t) &= P_1^2 \int _{\frac{1}{t^2}}^{1}  e^{-r^2t} \frac{\sin^2 (t\sqrt{r^2 + m^2 \log(1+r)} )}{r^2 + m^2 \log (1+r)}  dr  \leq P_1^2 \int _{\frac{1}{t^2}}^{1}   \frac{1 }{r^2 + m^2 \log (1+r)}  dr \nonumber \\
&\leq \frac{2}{m^2} P_1^2 \int _{\frac{1}{t^2}}^{1}   \frac{1 }{r}  dr = \frac{2}{m^2} P_1^2 (1 - \log t^{-2}) = \frac{2}{m^2}P_1^2 (1+ 2 \log t) \nonumber \\
 &\leq \frac{C}{m^2} P_1^2 \log t , \quad t \gg 1. \label{A3}
\end{align}
\noindent
\eqref{A1}, \eqref{A2} and \eqref{A3} imply the desired statement 
\begin{align}
\int _{|\xi|\leq 1} |\varphi(t,\xi)|^2 d\xi \leq \frac{C}{m^{2}}P_1^2 \log t, \quad t\gg 1
\end{align}
with some constant $C > 0$.

One next derives the lower bound estimates of \eqref{A1} in the case when $2\theta = 1$ and $n = 1$.\\
One starts with the following inequality:
\begin{align}\label{A5_03maio}
\int _{|\xi| \leq 1} & |\varphi(t,\xi)|^2 d\xi = P_1^2 \int _{0}^{1}  e^{-r^2t} \frac{\sin^2 (t\sqrt{r^2 + m^2 \log(1+r)} )}{r^2 + m^2 \log (1+r)}  dr \nonumber \\ 
&\geq P_1^2 \int _{0}^{1}  e^{-t(r^2 +m^2 \log(1+r)) } \frac{\sin^2 (t\sqrt{r^2 + m^2 \log(1+r)} )}{r^2 + m^2 \log (1+r)}  dr \nonumber \\ 
&= 2 P_1^2 \int _{0}^{1}   \frac{ e^{-t(r^2 +m^2 \log(1+r)) } \sin^2 (t\sqrt{r^2 + m^2 \log(1+r)} )}{\sqrt{t} \sqrt{r^2 + m^2 \log (1+r)} \left ( \frac{m^2}{r+1} +2r \right ) }  \frac{\sqrt{t} \left ( \frac{m^2}{r+1} +2r \right ) }{2\sqrt{r^2 + m^2 \log (1+r)}} dr.
\end{align}
%We observe that 
%$$ \frac{m^2}{r+1} +2r \approx 1 , \quad 0\leq r\leq 1.$$ 
Set $ \alpha = \sqrt{1 + m^2 \log 2}  $.  We apply the change of variable $ w= \sqrt{t} \sqrt{r^2 + m^2 \log(1+r)}  $ in \eqref{A5_03maio}. Then, one has% we have $$ d w =  \frac{ \sqrt{t} \left ( \frac{m^2}{r+1} +2r \right ) }{2 \sqrt{r^2 + m^2 \log(1+r)} } dr. $$
\begin{align}\label{A7_03maio}
\int _{|\xi| \leq 1} |\varphi(t,\xi)|^2 d\xi &\geq C^{*} P_1^2 \int _0^{\alpha\sqrt{t}} \frac{e^{-w^2 } \sin^2(\sqrt{t} w)}{w} dw, 
\end{align}
where one has used the estimate:
\[2\left(\frac{m^{2}}{r+1}+2r\right)^{-1} \geq 2(m^{2}+2)^{-1} =: C^{*}\quad (0 \leq r \leq 1).\]

Now, for $ t \geq 1 $, we observe  that
\begin{align}
\int _{\alpha \sqrt{t} }^\infty \frac{e^{-w^2 } \sin^2(\sqrt{t} w)}{w} dw &\leq \frac{1}{\alpha } t^{-\frac{1}{2}} \int _{\alpha \sqrt{t} }^\infty e^{-w^2 } dw \nonumber \\
&\leq  t^{-\frac{1}{2}} \int _{0 }^\infty e^{-w^2 } dw \nonumber \\
&= \frac{\sqrt{\pi}}{2} t^{-\frac{1}{2}}  .
\end{align}
Thus, applying Lemma \ref{LemmaA1} and the estimate just above, one obtains
\begin{align}
\int _0^{\alpha\sqrt{t}} \frac{e^{-w^2 } \sin^2(\sqrt{t} w)}{w} dw & = \int _0^{\infty} \frac{e^{-w^2 } \sin^2(\sqrt{t} w)}{w} dw - \int _{\alpha \sqrt{t} }^\infty \frac{e^{-w^2 } \sin^2(\sqrt{t} w)}{w} dw\nonumber \\
&\geq C\log t - \frac{\sqrt{\pi}}{2} t^{-\frac{1}{2}}, \label{A6_03maio}
%&= \log (t) \Big ( C - \frac{\sqrt{\pi}}{2}  \frac{t^{-\frac{1}{2}}}{\log(t)}  \Big ) . 
\end{align}
%Since $$ \lim _{t \rightarrow \infty } C - \frac{\sqrt{\pi}}{2}  \frac{t^{-\frac{1}{2}}}{\log(t)} = C, $$
%there exists $T>0$ such that 
%$$C - \frac{\sqrt{\pi}}{2}  \frac{t^{-\frac{1}{2}}}{\log(t)} \geq \frac{C}{2}, \quad t \geq T . $$
which implies %\eqref{A6_03maio}, we have 
\begin{align}\label{A9_maio03}
\int _0^{\alpha\sqrt{t}} \frac{e^{-w^2 } \sin^2(\sqrt{t} w)}{w} dw  \geq \frac{C}{2} \log t, \quad t \gg 1. 
\end{align}
\eqref{A7_03maio} and \eqref{A9_maio03} imply the desired estimate from below with a constant $C>0$:
$$ \int _{|\xi| \leq 1} |\varphi(t,\xi)|^2 d\xi \geq \frac{CC^{*}}{2} P_1^2  \log t, \quad t \gg 1.  $$
%Therefore, the lemma is proved. 
\hfill $\Box$

\subsubsection*{(III)\,Estimates for $n=2$ and $\theta =1$.}
Finally one treats the critical case for $n = 2$ and $\theta = 1$ to observe the growth property of the solution itself. In this case, one can get the following statement.
\begin{lem}\label{Lema6.5}
Let $ n=2$ and $  \theta = 1  $. There exist positive constants $C_1, C_2$ such that 
$$ \frac{C_1}{2+m^{2}}P_1^2 \log t  \leq  \int _{|\xi|\leq 1} |\varphi(t,\xi)|^2 d\xi \leq \frac{C_2}{1+m^{2}}P_1^2 \log t , \quad t\gg 1.$$
%where
%\[C_{0} := \left(2+m^{2}\right)^{-1}.\]
\end{lem}
\textit{Proof. } 
For $t>0$, for $n = 2$ and $\theta = 1$ we observe that 
\begin{align}
\int _{|\xi|\leq 1} |\varphi(t,\xi)|^2 d\xi &= P_1^2 \int _{0}^{1}  e^{-r^2t} \frac{\sin^2 (t\sqrt{r^2 + m^2 \log(1+r^{2})} )}{r^2 + m^2 \log (1+r^{2})} r  dr \\  
&=  P_1^2 \int _{0}^{\frac{1}{t}}  e^{-r^2t} \frac{\sin^2 (t\sqrt{r^2 + m^2 \log(1+r^{2})} )}{r^2 + m^2 \log (1+r^{2})} r  dr \nonumber \\ &+ P_1^2 \int _{\frac{1}{t}}^{1}  e^{-r^2t} \frac{\sin^2 (t\sqrt{r^2 + m^2 \log(1+r^{2})} )}{r^2 + m^2 \log (1+r^{2})} r  dr \nonumber \\
&=: P_1^2 ( M_1(t) + M_2(t)). \label{A8_maio03}
\end{align}
First of all, one knows that $ |\sin a | \leq a  $ for all $a \geq 0$. Thus 
\begin{align}
M_1(t) &= \int _{0}^{\frac{1}{t}}  e^{-r^2t} \frac{\sin^2 (t\sqrt{r^2 + m^2 \log(1+r^{2})} )}{r^2 + m^2 \log (1+r^{2})} r  dr \nonumber \\
&\leq t^2 \int _{0}^{\frac{1}{t}}  e^{-r^2t} r dr \leq  t^2 \int _{0}^{\frac{1}{t}}   r dr \nonumber \\ 
&= \frac{1}{2}.
\end{align}
In order to estimate $M_2(t)$, we observe that 
$$ \frac{r}{r^2 + m^2 \log (1+r^{2})} \approx \frac{1}{1+m^{2}}\frac{1}{r}, \quad 0 < r \leq 1 . $$
Then, the same steps to obtain estimates for $K_2(t)$  in \eqref{A3} can be applied to get 
\begin{align*}
M_2(t) \leq \frac{C}{1+m^{2}}\log t, \quad t \gg 1. 
\end{align*}
Applying the estimates to $M_1(t)$ and $M_2(t)$ in \eqref{A8_maio03}, we can prove the upper estimate
\begin{align}
\int _{|\xi|\leq 1} |\varphi(t,\xi)|^2 d\xi &\leq \frac{C}{1+m^{2}}P_1^2 \log (t), \quad t \gg 1.
\end{align}

To obtain the lower bound estimate we set $ \alpha = \sqrt{1+m^2 \log 2} $. Then, similarly to \eqref{A7_03maio} one can estimate as follows:
\begin{align}
\int _{|\xi|\leq 1} &|\varphi(t,\xi)|^2 d\xi = P_1^2 \int _{0}^{1}  e^{-r^2t} \frac{\sin^2 (t\sqrt{r^2 + m^2 \log(1+r^{2})} )}{r^2 + m^2 \log (1+r^{2})} r  dr \nonumber \\
&\geq P_1^2 \int _{0}^{1}   \frac{ e^{-t \Big ( r^2 +m^2 \log (1+r^2) \Big) }\sin^2 (t\sqrt{r^2 + m^2 \log(1+r^{2})} )}{ \sqrt{t}\sqrt{r^2 + m^2 \log(1+r^{2})} \Big ( 2  + \frac{m^2}{1+r^2} \Big ) }  \frac{ \sqrt{t} r \Big ( 2  + \frac{m^2}{1+r^2} \Big )  }{ \sqrt{r^2 + m^2 \log(1+r^{2})}} dr \nonumber \\
&\geq C_{0}P_1^2 \int _{0}^{\alpha \sqrt{t}} \frac{e^{-w^2} \sin (\sqrt{t} w )}{w } dw, 
\end{align}
where we have just used the change of variable $w= \sqrt{t}\sqrt{r^2+m^2 \log (1+r^2) }$ and the fact  
$$\left( 2 + \frac{m^2}{1+r^2} \right)^{-1} \geq (2+m^{2})^{-1} =: C_{0}\quad (0 \leq r \leq 1).$$

As we have already proved in \eqref{A9_maio03}, we have the desired estimate:
$$\int _{|\xi|\leq 1} |\varphi(t,\xi)|^2 d\xi \geq C_{0}P_1^2 \log t, \quad t \gg 1.  $$
\hfill
$\Box$

Now, let us check the validity of the statement of Theorem 1.1.\\
{\it Proof of Theorem 1.1.}
The proof of Theorem 1.1 is a direct consequence of \eqref{solu-for}, \eqref{est-F1c}, \eqref{est-F2b}, \eqref{est-F3-A3}, \eqref{mf-10},  \eqref{sub-hi},
\eqref{hyper-hi} and \eqref{perfil-high}.
\hfill
$\Box$
%\section{ Main  results}
%Based on the previous estimates  we can prove the following theorems.
%\begin{theo}\label{Prop-1}
%Let $0\leq \theta \leq 1$ and  $n \geq 1$,  $\varphi$  the function defined in \eqref{Perfil}. Assuming that $ u_1 \in L^{1}(\R^n) \cap L^2(\R^n)$, there exists a positive constant $C=C_{n,\delta_0,\theta}$ such that 
%\begin{align*}
%\| u(t, \cdot) - \mathcal{F}^{-1}_{\xi \to x } (\varphi (t, \xi) ) (\cdot)  \|^2 \leq C \| u_1 \| _{1, \theta}^2 t^{-\frac{n}{2}}, \quad t \gg 1. 
%\end{align*}
%\end{theo}
%\textit{Proof.} 

Therefore, based on Theorem 1.1, the Plancherel Theorem and the standard inequality
\begin{equation}\label{triangle}
\| \varphi(t,\cdot)\|  - \|\hat{u}(t,\cdot) - \varphi(t,\cdot)\| \leq \|\hat{u}(t,\cdot)\| \leq   \|\hat{u}(t,\cdot) - \varphi(t,\cdot)\| + \| \varphi(t,\cdot)\|,\;\; t\geq 1,
\end{equation}
one can prove Theorem 1.2.\\
{\it Proof of theorem 1.2.}  
The proof of item \rm{(i)} follows by using \eqref{triangle}, Theorem 1.1, \eqref{perfil-above}, \eqref{Lem-perfil-low}. Similarly, the proof of items \rm{(ii)} and \rm{(iii)} follows from Theorem 1.1 combined with Lemmas \ref{Lema6.3}, \ref{Lema6.4} and \ref{Lema6.5}. 
\hfill $\Box$

\par
\vspace{0.5cm}
\noindent{\em Acknowledgements.}
\smallskip
The work of the first  author (A. Piske) was financed in part by the Coordena\c{c}\~ao de Aperfei\c{c}oamento de Pessoal de N\'ivel Superior - Brasil (CAPES) - Finance Code 001. The work of the second author (R. C. Char\~ao) was partially supported by PRINT/CAPES-Process 88881.310536/2018-00.  The work of the third author (R. Ikehata) was supported in part by Grant-in-Aid for Scientific Research (C) 20K03682 of JSPS.

%%%%%%%%%%%%%%%%%%%%%%%%%%%%%%%%%%%%%%%%%%%%%%%%%%%%%%%%%%%%%%%%%%%%%%%%%%%%%%%%%%%%%%%%%

%%%%%%%%%%%%%%%%%%%%%%%%%%%%%%%%%%%%%%%%%%%

\end{document}